\hsize=6 truein
\nopagenumbers
\baselineskip=18pt plus 5pt minus 4pt
\hoffset= .15truein

\font\tenmsya=msam10
\font\sevenmsya=msam7
\font\fivemsya=msam5
\newfam\msyafam 
\textfont\msyafam=\tenmsya
\scriptfont\msyafam=\sevenmsya
\scriptscriptfont\msyafam=\fivemsya

\font\tenmsyb=msbm10
\font\sevenmsyb=msbm7
\font\fivemsyb=msbm5
\newfam\msybfam \def\msyb{\fam\msybfam\tenmsyb}
\textfont\msybfam=\tenmsyb \scriptfont\msybfam=\sevenmsyb
\scriptscriptfont\msybfam=\fivemsyb

\def\magstep#1{\ifcase#1 \@m\or 1200\or 1440\or 1728\or 2074\or 2488
       \or 2986\or 3583\or 4300\or 5160\fi\relax}

\def\Itemitem{\hangindent2\parindent\textindent}

   \def\calH{\ifmmode{\cal H}\else$\cal H$\fi}
   \def\calO{\ifmmode{\cal O}\else$\cal O$\fi}
   \def\calP{\ifmmode{\cal P}\else$\cal P$\fi}
   \def\calS{\ifmmode{\cal S}\else$\cal S$\fi}
   \def\calT{\ifmmode{\cal T}\else$\cal T$\fi}
   \def\calM{\ifmmode{\cal M}\else$\cal M$\fi}
   \def\calN{\ifmmode{\cal N}\else$\cal N$\fi}
  \def\n{\ifmmode{\msyb N}\else$\msyb N$\fi}
  \def\r{\ifmmode{\msyb R}\else$\msyb R$\fi}
  \def\z{\ifmmode{\msyb Z}\else$\msyb Z$\fi}
  \def\q{\ifmmode{\msyb Q}\else$\msyb Q$\fi}
\def\c{\ifmmode{\msyb C}\else$\msyb C$\fi}

\hfuzz=.5in

\def\frac#1#2{{#1\over #2}}


\centerline{\bf Classification of 3-dimensional isolated rational
hypersurface singularities with $\c^*$-action}


\bigskip

\hskip44pt {\bf Stephen S.-T. Yau$^*$}
\vskip-6pt
\hskip00pt {\bf Department of Mathematics, Statistics}
\vskip-6pt
\hskip10pt {\bf and Computer Science (M/C 249)}
\vskip-6pt
\hskip13pt {\bf University of Illinois at Chicago}
\vskip-6pt
\hskip30pt {\bf 851 South Morgan Street}
\vskip-6pt
\hskip16pt {\bf Chicago, IL 60607--7045, U.S.A.}
\vskip-6pt
\hskip40pt {\bf e--mail: yau@uic.edu}

\vskip-90pt

\hskip266pt {\bf Yung Yu$^{**}$}
\vskip-6pt
\hskip226pt {\bf Department of Mathematics}
\vskip-6pt
\hskip214pt {\bf National Cheng Kung University}
\vskip-6pt
\hskip237pt {\bf Tainan, Taiwan, R.O.C.}
\vskip-6pt
\hskip212pt {\bf e--mail: yungyu@mail.ncku.edu.tw}

\vskip7in

\noindent $^*$\hskip8pt Research partially supported by NSF, NSA, U.S.A.

\noindent $^{**}$ Research partially supported by NSC, R.O.C.

\vfill\eject

\footline={\hss\tenrm\folio\hss}
\pageno=1

\noindent {\bf 1. Introduction}

\smallskip

In [Art] Artin first introduced the definition of rational surface
singularity.
He classified all rational surface singularities embeddable in $\c^3$.
These are precisely those Du Val Singularities in $\c^3$ defined by one
of the following polynomial equations:

\medskip
\itemitem{}
$A_n$:\quad $x^2+y^2+z^{n+1}$, for $n\geq 1$
\smallskip

\itemitem{}
$D_n$:\quad $x^2+y^2z+z^{n-1}$, for $n\geq 4$
\smallskip
\itemitem{}
$E_6$:\quad $x^2+y^3+z^4$
\smallskip
\itemitem{}
$E_7$:\quad $x^2+y^3+yz^3$
\smallskip
\itemitem{}
$E_8$:\quad $x^2+y^3+z^5$

\medskip

\noindent It is well known that any canonical singularity (i.e. singularity
that occurs
in a canonical model of a surface of general type)
is analytically isomorphic to one of the rational double points
listed above.

In [Bu] Burns defined higher dimensional rational singularity
as follows.
Let $(V,p)$ be a $n$-dimensional isolated singularity.
Let $\pi:M\to V$ be a resolution of singularity.
$p$ is said be a rational singularity if $R^i\pi_*\calO_M=0$
for $1\leq i\leq n-1$.
In [Ya4], Yau shows for Gorenstein singularities that it is sufficient to
require $R^{n-1}\pi_*\calO_M=0$.
He [Ya2] proves that
$R^{n-1}\pi_*\calO_M\cong H^0(V-\{p\},\Omega^n)/L^2(V-\{p\},\Omega^n)$
where $\Omega^n$ is the sheaf of germs of holomorphic
$n$-forms and $L^2(V-\{p\},\Omega^n)$ is the space of holomorphic
$n$ forms on $V-\{p\}$ which are $L^2$-integrable.
The geometric genus $p_g$ of the singularity $(V,p)$ is defined to be
$$
p_g:=\hbox{dim} R^{n-1}\pi_*\calO_M=\hbox{dim} H^0(V-\{p\},\Omega^n)
/L^2(V-\{p\},\Omega^n)
$$
It turns out that $p_g$ is an important invariant of $(V,p)$.

In [Ya-Yu], we give algebraic classification of rational CR structures
on the topological 5-sphere with transversal holomorphic $S^1$-action
in $\c^4$.
Here, algebraic classification of compact strongly pseudoconvex CR
manifolds $X$ means classification up to algebraic equivalence,
i.e. roughly up to isomorphism of the normalization of the complex
analytic variety $V$ which has $X$ as boundary.
The problem is intimately related to the study of 3-dimensional
isolated rational weighted homogeneous hypersurface singularities with link
homeomorphic to $S^5$.
For this, we need the classification of 3-dimensional isolated rational
hypersurface singularities with a $\c^*$-action.
This list is only available at the homepage of one of us.
Since there is a desire for a complete list of this classification
(cf. Theorem 3.3), we decide to publish it for the convenience of readers.

The idea of our proof is very easy.
If $h(z_0,z_1,z_2,z_3)$ is a weighted homogeneous polynomial in
$\c^4$ and $V=\{z\in\c^4:h(z)=0\}$ has an isolated singularity at the
origin, then Kouchnirenko [Ko] and Orlik-Randell [Or-Ra] observed that
$V$ can be deformed into one of the nineteen classes of weighted homogeneous
singularities listed in \S 2 while keeping the differential structure of
the link $K_V:=S^7\cap V$ constant.
We prove that the above deformation is actually a
deformation that preserves weights and embedded
topological type without changing weights.
By a theorem of Merle-Teissier [Me-Te], the geometric genus $p_g$ of the
singularity can be expressed in terms of its weights.
The MAPLE program [Ch] helps us to finish the classification.
In fact if we use the similar method as above, we can also classify the
weighted homogenous rational surface singularities embeddable in $\c^3$,
which are exactly $A_n,D_n,E_6,E_7$, and $E_8$ singularities described above.

In \S2, we shall give a classification (up to deformation which
preserves weights) of weighted homogeneous polynomials of 4 variables
with isolated singularity at the origin.
This list was obtained first by Kouchnirenko [Ko] and
Orlik-Randell [Or-Ra] (see also [Ka]) independently.
In \S3, we classify all 3-dimensional isolated rational hypersurface
singularities with $\c^*$-action.

\bigskip
\bigskip

\noindent {\bf 2. Classification of weighted homogeneous polynomials
in four variables with isolated singularity at the origin}

\smallskip

Orlik and Wagreich [Or-Wa] and Arnold [Ar] showed that if
$h(z_0,z_1,z_2)$ is a weighted homogeneous polynomial in $\c^3$ and
$V=\{z\in\c^3:h(z)=0\}$ has an isolated singularity at the origin, then
$V$ can be deformed into one of the following seven classes of weighted
homogeneous singularities while keeping the differential structure
of the link $K_V:=S^5\cap V$ constant.
Let $(w_0,w_1,w_2)=(wt(z_0),wt(z_1),wt(z_2))$ be the weight type and $\mu$ be
the Milnor number.

\medskip

\noindent {\bf Class I}\quad $\{z^{a}_0+z^{b}_1+z^{c}_2=0\}$,
$(w_0,w_1,w_2)=(a,b,c)$, $\mu=(a-1)(b-1)(c-1)$.

\medskip

\noindent {\bf Class II}\quad $\{z^{a}_0+z^{b}_1+z_1z^{c}_2=0\}$,
$(w_0,w_1,w_2)=(a,b,\frac{bc}{b-1})$, $\mu=(a-1)(bc-b+1)$.

\medskip

\noindent {\bf Class III}\quad $\{z^{a}_0+z^{b}_1z_2+z_1z^{c}_2=0\}$,
$(w_0,w_1,w_2)=(a,\frac{bc-1}{c-1},\frac{bc-1}{b-1})$,
$\mu=(a-1)bc$.

\medskip

\noindent {\bf Class IV}\quad $\{z^{a}_0+z^{b}_1z_2+z_0z^{c}_2=0\}$,
$(w_0,w_1,w_2)=(a,\frac{abc}{ac-a+1},\frac{ac}{a-1})$,
$\mu=ac(b-1)+a-1$.

\medskip

\noindent {\bf Class V}\quad $\{z^{a}_0z_1+z^{b}_1z_2+z_0z^{c}_2=0\}$,
$(w_0,w_1,w_2)=(\frac{abc+1}{bc-c+1},
        \frac{abc+1}{ac-a+1},
        \frac{abc+1}{ab-b+1})$,

\hskip30pt $\mu=abc$.

\medskip

\noindent {\bf Class VI}\quad
$\{z^{a}_0+z_0z^{b}_1+z_0z^{c}_2+z_1^{b_1}z_2^{b_2}=0\}$,
where $(a-1)(bb_2+cb_1)=abc$,

\hskip34pt $(w_0,w_1,w_2)=(a,\frac{ab}{a-1},\frac{ac}{a-1})$,
$\mu=\frac{(ab-a+1)(ac-a+1)}{a-1}$.

\medskip

\noindent {\bf Class VII}\quad
$\{z^{a}_0z_1+z_0z^{b}_1+z_0z^{c}_2+z_1^{b_1}z_2^{b_2}=0\}$,
where $(a-1)(bb_2+cb_1)=c(ab-1)$,

\hskip34pt $(w_0,w_1,w_2)=(\frac{ab-1}{b-1},\frac{ab-1}{a-1},
\frac{c(ab-1)}{b(a-1)})$,
$\mu=\frac{a(abc-ab+b-c)}{a-1}$.

Recall that two isolated hypersurface singularities
$(V,0), (W,0)$ in $\c^{n+1}$ are said to have the same topological type
if $(\c^{n+1},V,0)$ is homeomorphic to $(\c^{n+1},W,0)$ (cf. [Ya1]).

In [Xu-Ya1], we prove that the above deformation is actually a topological
type constant deformation without changing weights.
Therefore any weighted homogeneous singularity has the same topological type
of one of the seven classes above.

If $h(z_0,z_1,z_2,z_3)$ is a weighted homogeneous polynomial in $\c^4$ and
$V=\{z\in\c^4:h(z)=0\}$ has an isolated singularity at the origin, then
Kouchnirenko [Ko] and Orlik and Randell [Or-Ra] observed that $V$ can be
deformed into one of the following nineteen classes of weighted homogeneous
singularities below while keeping the differential structure of the link
$K_V:=S^7\cap V$ constant (the meaning of the linear forms $\alpha$ in the list
will be explained later).
Let $w_i=wt(z_i)$ and $\mu$ be the Milnor number.

\medskip

\noindent {\bf Type I}\quad $\{z^{a}_0+z^{b}_1+z^{c}_2+z^{d}_3=0\}$,
$\frac{x}{a}+\frac{y}{b}+\frac{z}{c}+\frac{w}{d}=\alpha(x,y,z,w)$,
$(w_0,w_1,w_2,w_3)=(a,b,c,d)$,

\hskip26pt $\mu=(a-1)(b-1)(c-1)(d-1)$.

\medskip

\noindent {\bf Type II}\quad
$\{z^{a}_0+z^{b}_1+z^{c}_2+z_2z^{d}_3=0\}$,
$\frac{x}{a}+\frac{y}{b}+\frac{z}{c}+\frac{(c-1)w}{cd}=\alpha(x,y,z,w)$,
$(w_0,w_1,w_2,w_3)=(a,b,c,\frac{cd}{c-1})$,

\hskip32pt $\mu=(a-1)(b-1)[c(d-1)+1]$.

\medskip

\noindent {\bf Type III}\quad
$\{z^{a}_0+z^{b}_1+z^{d}_2z_3+z_2z^{a_4}_3=0\}$,
$\frac{x}{a}+\frac{y}{b}+\frac{(d-1)z}{cd-1}+\frac{(c-1)w}{cd-1}=\alpha(x,y,z,w)$,

\hskip38pt $(w_0,w_1,w_2,w_3)=(a,b,\frac{cd-1}{d-1},\frac{cd-1}{c-1})$,
$\mu=(a-1)(b-1)cd$.

\medskip

\noindent {\bf Type IV}\quad
$\{z^{a}_0+z_0z^{b}_1+z_2^{c}+z_2z^d_3=0\}$,
$\frac{x}{a}+\frac{(a-1)y}{ab}+\frac{z}{c}+\frac{(c-1)w}{cd}=\alpha(x,y,z,w)$,

\hskip38pt $(w_0,w_1,w_2,w_3)=(a,\frac{ab}{a-1},c,\frac{cd}{c-1})$,
$\mu=[a(b-1)+1][c(d-1)+1]$.

\medskip

\noindent {\bf Type V}\quad
$\{z^{a}_0z_1+z_0z^{b}_1+z^{c}_2+z_2z_3^{d}=0\}$,
$\frac{(b-1)x}{ab-1}+\frac{(a-1)y}{ab-1}+\frac{z}{c}+\frac{(c-1)w}{cd}=\alpha(x,y,z,w)$,

\hskip32pt $(w_0,w_1,w_2,w_3)=(\frac{ab-1}{b-1},
        \frac{ab-1}{a-1},c,
        \frac{cd}{c-1})$,
$\mu=ab[c(d-1)+1]$.

\medskip

\noindent {\bf Type VI}\quad
$\{z^{a}_0z_1+z_0z^{b}_1+z^{c}_2z_3+z_2z_3^{d}=0\}$,
$\frac{(b-1)x}{ab-1}+\frac{(a-1)y}{ab-1}+\frac{(d-1)z}{cd-1}+\frac{(c-1)w}{dc}=\alpha(x,y,z,w)$,

\hskip38pt $(w_0,w_1,w_2,w_3)=(\frac{ab-1}{b-1},
        \frac{ab-1}{a-1},
        \frac{cd-1}{d-1},
        \frac{cd-1}{c-1})$,
$\mu=abcd$.

\medskip
\noindent {\bf Type VII}\quad
$\{z^{a}_0+z^{b}_1+z_1z^{c}_2+z_2z_3^{d}=0\}$,
$\frac{x}{a}+\frac{y}{b}+\frac{(b-1)z}{bc}+\frac{[b(c-1)+1]w}{bcd}=\alpha(x,y,z,w)$,

\hskip42pt $(w_0,w_1,w_2,w_3)=(a,b,\frac{bc}{b-1},
        \frac{bcd}{b(c-1)+1})$,
$\mu=(a-1)[bc(d-1)+b-1]$.

\medskip
\noindent {\bf Type VIII}\quad
$\{z^{a}_0+z^{b}_1+z_1z^{c}_2+z_1z_3^{d}+z_2^pz_3^q=0$,
$\frac{p(b-1)}{bc}+\frac{q(b-1)}{bc}=1\}$,

\hskip46pt $\frac{x}{a}+\frac{y}{b}+\frac{(b-1)z}{bc}+\frac{(b-1)w}{bd}=\alpha(x,y,z,w)$,

\hskip46pt $(w_0,w_1,w_2,w_3)=(a,b,\frac{bc}{b-1},
        \frac{bd}{b-1})$,
$\mu=\frac{(a-1)[b(c-1)+1][b(d-1)+1]}{b-1}$.

\medskip
\noindent {\bf Type IX}\quad
$\{z^{a}_0+z^{b}_1z_3+z^{c}_2z_3+z_1z_3^{d}+z^p_1z^q_2=0$,

\hskip38pt $\frac{p(d-1)}{bd-1}+\frac{qb(d-1)}{c(bd-1)}=1\}$,
$\frac{x}{a}+\frac{(d-1)y}{bd-1}+\frac{b(d-1)z}{c(bd-1)}+\frac{(b-1)w}{bd-1}=\alpha(x,y,z,w)$,

\hskip38pt $(w_0,w_1,w_2,w_3)=(a,\frac{bd-1}{d-1},
        \frac{c(bd-1)}{b(d-1)},
        \frac{bd-1}{b-1})$,
$\mu=\frac{(a-1)d[c(bd-1)-b(d-1)]}{d-1}$.

\medskip
\noindent {\bf Type X}\quad
$\{z^{a}_0+z^{b}_1z_2+z^{c}_2z_3+z_1z_3^{d}=0\}$,
$\frac{x}{a}+\frac{[d(c-1)+1]y}{bcd+1}+\frac{[b(d-1)+1]z}{bcd+1}
+\frac{[c(b-1)+1]w}{bcd+1}=\alpha(x,y,z,w)$,

\hskip32pt $(w_0,w_1,w_2,w_3)=(a,\frac{bcd+1}{d(c-1)+1},
               \frac{bcd+1}{b(d-1)+1},
              \frac{bcd+1}{c(b-1)+1})$,
$\mu=(a-1)bcd$.

\medskip
\noindent {\bf Type XI}\quad
$\{z^{a}_0+z_0z^{b}_1+z_1z^{c}_2+z_2z_3^{d}=0\}$,
$\frac{x}{a}+\frac{(a-1)y}{ab}+\frac{[a(b-1)+1]z}{abc}
+\frac{[ab(c-1)+(a-1)]w}{abcd}=\alpha(x,y,z,w)$,

\hskip36pt $(w_0,w_1,w_2,w_3)=(a,\frac{ab}{a-1},
        \frac{abc}{a(b-1)+1},
        \frac{abcd}{ab(c-1)+(a-1)})$,
$\mu=abc(d-1)+a(b-1)+1$.

\medskip
\noindent {\bf Type XII}\quad
$\{z^{a}_0+z_0z^{b}_1+z_0z^{c}_2+z_1z_3^{d}+z^p_1z^q_2=0$,
$\frac{p(a-1)}{ab}+\frac{q(a-1)}{ac}=1\}$,

\hskip40pt $\frac{x}{a}+\frac{(a-1)y}{ab}+\frac{(a-1)z}{ac}+\frac{[a(b-1)+1]w}{abd}=\alpha(x,y,z,w)$,

\hskip40pt $(w_0,w_1,w_2,w_3)=(a,\frac{ab}{a-1},
        \frac{ac}{a-1},
        \frac{abd}{a(b-1)+1})$,
$\mu=\frac{(a(c-1)+1)(ab(d-1)+a-1)}{a-1}$.

\medskip
\noindent {\bf Type XIII}\quad
$\{z^{a}_0+z_0z^{b}_1+z_1z^{c}_2+z_1z_3^{d}+z^p_2z^q_3=0$,
$\frac{p[a(b-1)+1]}{abc}+\frac{q[a(b-1)+1]}{abd}=1\}$,

\hskip42pt $\frac{x}{a}+\frac{(a-1)y}{ab}+\frac{[a(b-1)+1]z}{abc}
+\frac{[a(b-1)+1]w}{abd}=\alpha(x,y,z,w)$,

\hskip42pt $(w_0,w_1,w_2,w_3)=(a,\frac{ab}{a-1},
        \frac{abb}{a(b-1)+1},
        \frac{abd}{a(b-1)+1})$,
$\mu=\frac{[ab(c-1)+a-1][ab(d-1)+a-1]}{a(b-1)+1}$.

\medskip
\noindent {\bf Type XIV}\quad
$\{z^{a}_0+z_0z^{b}_1+z_0z^{c}_2+z_0z_3^{d}+z^p_1z^q_2+z_2^rz_3^s=0$,
$\frac{p(a-1)}{ab}+\frac{q(a-1)}{ac}=1$ and

\hskip46pt $\frac{r(a-1)}{ac}+\frac{s(a-1)}{ad}=1\}$,
$\frac{x}{a}+\frac{(a-1)y}{ab}+\frac{(a-1)z}{ac}+\frac{(a-1)w}{ad}=\alpha(x,y,z,w)$,

\hskip46pt $(w_0,w_1,w_2,w_3)=(a,\frac{ab}{a-1},
               \frac{ac}{a-1},\frac{ad}{a-1})$,
$\mu=\frac{[a(b-1)+1][a(c-1)+1][a(d-1)+1]}{(a-1)^2}$.

\medskip
\noindent {\bf Type XV}\quad
$\{z^{a}_0z_1+z_0z^{b}_1+z_0z^{c}_2+z_2z_3^{d}+z^p_1z^q_2=0$,
$\frac{p(a-1)}{ab-1}+\frac{qb(a-1)}{c(ab-1)}=1\}$,

\hskip40pt $\frac{(b-1)x}{ab-1}+\frac{(a-1)y}{ab-1}+\frac{b(a-1)z}{c(ab-1)}
+\frac{[c(ab-1)-b(a-1)]w}{cd(ab-1)}=\alpha(x,y,z,w)$,

\hskip40pt $(w_0,w_1,w_2,w_3)=(\frac{ab-1}{b-1},\frac{ab-1}{a-1},
        \frac{c(ab-1)}{b(a-1)},
        \frac{cd(ab-1)}{c(ab-1)-b(a-1)})$,
$\mu=\frac{a[c(d-1)(ab-1)+b(a-1)]}{a-1}$.

\medskip
\noindent {\bf Type XVI}\quad
$\{z^{a}_0z_1+z_0z^{b}_1+z_0z^{c}_2+z_0z_3^{d}+z^p_1z^q_2+z^r_2z^s_3=0$,
$\frac{p(a-1)}{ab-1}+\frac{qb(a-1)}{c(ab-1)}=1$ and

\hskip44pt $\frac{rb(a-1)}{c(ab-1)}+\frac{sb(a-1)}{d(ab-1)}=1\}$,
$\frac{(b-1)x}{ab-1}+\frac{(a-1)y}{ab-1}+\frac{b(a-1)z}{c(ab-1)}
+\frac{b(a-1)w}{d(ab-1)}=\alpha(x,y,z,w)$,

\hskip44pt $(w_0,w_1,w_2,w_3)=(\frac{ab-1}{b-1},\frac{ab-1}{a-1},
            \frac{c(ab-1)}{b(a-1)},
            \frac{d(ab-1)}{b(a-1)})$,
$\mu=\frac{a[c(ab-1)-b(a-1)]
   [d(ab-1)-b(a-1)]}{b(a-1)^2}$.

\medskip
\noindent {\bf Type XVII}\quad
$\{z^{a}_0z_1+z_0z^{b}_1+z_1z^{c}_2+z_0z_3^{d}+z^p_1z^q_3+z^r_0z^s_2=0$,
$\frac{p(a-1)}{ab-1}+\frac{qb(a-1)}{d(ab-1)}=1$ and

\hskip48pt $\frac{r(b-1)}{ab-1}+\frac{sa(b-1)}{c(ab-1)}=1\}$,
$\frac{(b-1)x}{ab-1}+\frac{(a-1)y}{ab-1}+\frac{a(b-1)z}{c(ab-1)}
+\frac{b(a-1)w}{d(ab-1)}=\alpha(x,y,z,w)$,

\hskip48pt $(w_0,w_1,w_2,w_3)=(\frac{ab-1}{b-1},
        \frac{ab-1}{a-1},
        \frac{c(ab-1)}{a(b-1)},
        \frac{d(ab-1)}{c(a-1)})$,
$\mu=\frac{[c(ab-1)-a(b-1)][d(ab-1)-b(a-1)]}
      {(a-1)(b-1)}$.

\medskip
\noindent {\bf Type XVIII}\quad
$\{z^{a}_0z_2+z_0z^{b}_1+z_1z^{c}_2+z_1z_3^{d}+z^p_2z^q_3=0$,
$\frac{p[a(b-1)+1]}{abc+1}+
\frac{qc[a(b-1)+1]}{d(abc+1)}=1\}$,

\hskip54pt $\frac{[b(c-1)+1]x}{abc+1}+\frac{[c(a-1)+1]y}{abc+1}
+\frac{[a(b-1)+1]z}{c(abc+1)}+\frac{c[a(b-1)+1]w}{d(abc+1)}=\alpha(x,y,z,w)$,

\hskip54pt $(w_0,w_1,w_2,w_3)=(\frac{abc+1}{b(c-1)+1},
        \frac{abc+1}{c(a-1)+1},
        \frac{abc+1}{a(b-1)+1},
        \frac{d(abc+1)}{c[a(b-1)+1]})$,
$\mu=\frac{ab[abc(d-1)+c(a-1)+d]}{a(b-1)+1}$.

\medskip
\noindent {\bf Type XIX}\quad
$\{z^{a}_0z_3+z_0z^{b}_1+z^{c}_2z_3+z_2z_3^{d}=0$,

\hskip44pt $\frac{[b(d(c-1)+1)-1]x}{abcd-1}+\frac{[d(c(a-1)+1)-1]y}{abcd-1}
+\frac{[a(b(d-1)+1)-1]z}{abcd-1}+\frac{[c(a(b-1)+1)-1]w}{abcd-1}=\alpha(x,y,z,w)$,

\hskip44pt $(w_0,w_1,w_2,w_3)=(\frac{abcd-1}{b[d(c-1)+1]-1},
        \frac{abcd-1}{d[c(a-1)+1]-1},
        \frac{abcd-1}{a[b(d-1)+1]-1},
        \frac{abcd-1}{c[a(b-1)+1]-1})$,

\hskip44pt $\mu=abcd$.

\medskip

\noindent {\bf Theorem 2.1.}\quad
Suppose $h(z_0,z_1,z_2,z_3)$ is a polynomial and
$V_k=\{(z_0,z_1,z_2,z_3)\in\c^4:h(z_0,z_1,z_2,z_3)=0\}$
has an isolated singularity at $0$.
Then $h(z_0,z_1,z_2,z_3)=f(z_0,z_1,z_2,z_3)+g(z_0,z_1,z_2,z_3)$
where $f$ is one of the 19 classes above with only isolated singularity at $0$
and $f$ and $g$ have no monomial in common.
In this case let
$V_f=\{(z_0,z_1,z_2,z_3)\in\c^4:f(z_0,z_1,z_2,z_3)=0\}$ and let
$$
K_f=V_f\cap S^7,\quad K_h=V_h\cap S^7.
$$
Then $K_f$ is equivariantly diffeomorphic to $K_h$.

\medskip
\noindent {\bf Proof.}\quad
If none of the monomials in
$\{z^{a_0}_0,z_0^{a_0}z_1,z_0^{a_0}z_2,z_0^{a_0}z_3\}$
appears in $h(z_0,z_1,z_2,z_3)$, then
$\frac{\partial h}{\partial z_j}(z_0,0,0,0)=0$, $0\leq j\leq 3$.
This contradicts the fact that $h$ has an isolated singularity at $0$.
Therefore, one of the monomials in
$\{z^{a_0}_0,z_0^{a_0}z_1,z_0^{a_0}z_2,z_0^{a_0}z_3\}$
appears in $h$.
Similarly one of the monomials in each of the following sets appears in
$h:\{z_0z_1^{a_1},z_1^{a_1},z_1^{a_1}z_2,z_1^{a_1}z_3\}$,
$\{z_0z_2^{a_2},z_1z_2^{a_2},z_2^{a_2},z_2^{a_2}z_3\}$,
$\{z_0z_3^{a_3},z_1z_3^{a_3},z_2z_3^{a_3},z_3^{a_3}\}$.
Taking a monomial from each of the 4 sets above, we get 256 polynomials.
One can check that these 256 polynomials are equivalent to one of the 19
classes above up to permutation of coordinates.
Notice that in Type VIII, for example, the monomial $z_2^pz_3^q$ is added to
make sure that $f$ has an isolated singularity at $0$.
Obviously if $h$ is weighted homogeneous of type $(w_0,w_1,w_2,w_3)$,
then so are $f$ and $g$.

The proof of Theorem 3.1.4 in [Or-Wa] shows that $K_f$ is equivariantly
diffeomorphic to $K_h$.
\hfill{Q.E.D.}

\medskip

We shall use the theory developed in [Xu-Ya1] and [Xu-Ya2] to show that
$(V_f,0)$ and $(V_h,0)$ have the same embedded topological type.

\medskip

\noindent {\bf Definition.}\quad
Given a real manifold $B$ of dimension $m$, and a family
$\{(M_t,N_t): t\in B_t$, $N_t$ is a closed submanifold of a compact
differentiable manifold $M_t\}$,
we say that $(M_t,N_t)$ depends $C^\infty$ on $t$ and that
$\{(M_t,N_t):t\in B\}$ is a $C^\infty$ family of compact manifolds with
submanifolds, if there is a $C^\infty$ mainifold $\calM$,
a closed submanifold $\calN$ and a $C^\infty$ map $w$ from $\calM$ onto $B$
such that $\overline{w}:=w|\calN$ is also a $C^\infty$ map from $\calN$
onto $B$ satisfying the following conditions
\itemitem{(i)}
$M_t=w^{-1}(t)\supseteq N_t=\overline{w}^{-1}(t)$
\itemitem{(ii)}
The rank of the Jacobian of $w$ (respectively $\overline{w}$)
is equal to $m$ at every point of $\calM$ (respectively $\calN$).

\medskip

\noindent {\bf Theorem 2.2.} (e.g., [Xu-Ya1])\quad
Let $((\calM,\calN), (w,\overline{w}))$ be a $C^\infty$ family of compact
manifolds with submanifolds, with $B$ connected.
Then $(M_t,N_t)=(w^{-1}(t),\overline{w}^{-1}(t))$ is diffeomorphic to
$(M_{t_0},N_{t_0})$ for any $t,t_0\in B$.

Now we fix weights $(w_0,\ldots,w_n)$ with $w_i\geq 2$.
Suppose that there is a weighted homogeneous polynomial
$f(z_0,\ldots,z_n)$ with the weights $(w_0,\ldots,w_n)$ such that $f$ has an
isolated singularity at the origin.
Let $\Delta$ be the intersection of the plane
$\sum\limits^n_{i=0}\frac{x_i}{w_i}=1$ with the first quadrant of $\r^{n+1}$.
Let $\c[\Delta]=\{f\in\c[z_0,\ldots,z_n]:$ supp $f\subset\Delta\}$
where supp $f=\{(d_0,\ldots,d_n)\in\r^{n+1}:
z^{d_0}_0z^{d_1}_1\cdots z^{d_n}_n$ occurs in $f\}$.
Let $N$ be the number of the integer points which are in $\Delta$.
There is a canonical correspondence between $\c[\Delta]$ and $\c^N$.
Thus we may introduce a Zariski topology on $\c[\Delta]$.

\medskip

\noindent {\bf Theorem 2.3.}\quad
Notation as above. Let
$$
U=\{f\in \c[\Delta]:f\,\,\hbox{has an isolated singularity at the origin}\,\}
$$
Then $U$ is a nonempty Zariski open set of $\c[\Delta]$.

The proof of the previous theorem as well as the following theorem is the
same as those of Theorem 3.4 and Theorem 3.5 in [Xu-Ya1] respectively.

\medskip

\noindent {\bf Theorem 2.4.}\quad
Suppose that $f(z_0,\ldots,z_n)$ and $g(z_0,\ldots,z_n)$ are weighted
homogeneous polynomials with the same weights $(w_0,w_1,\ldots,w_n)$.
If the variety $V$ of $f$ and the variety $W$ of $g$ have an isolated
singularity at the origin, then $(\c^{n+1},V,0)$ is homeomorphically equivalent
to $(\c^{n+1},W,0)$.

\medskip

\noindent {\bf Corollary 2.5.}\quad
Suppose that $h(z_0,z_1,z_2,z_3)$ is a weighted homogeneous polynomial with
weights $(w_0,w_1,w_2,w_3)$ and the variety $h^{-1}(0)$ has an isolated
singularity at the origin.
Then $h=f+g$ where $f$ and $g$ have no monomials in common, $f$ is one
of the nineteen classes above and $f$ and $g$ are weighted homogeneous
of type $(w_0,w_1,w_2,w_3)$.

Moreover $h^{-1}(0)$ and $f^{-1}(0)$ have the same embedded topological type.

\bigskip
\noindent {\bf 3. 3-dimensional isolated rational hypersurface singularities
with $\c^*$-action}

\medskip

\noindent {\bf Definition 3.1.}\quad
Let $(V,0)$ be an $n$-dimensional variety with isolated singularity at $0$.
The geometric genus $p_g(V,0)$ of the singularity is defined to be
dim $H^{n-1}(M,\calO)$ where $M$ is a resolution of the singularity
$(V,0)$.
$(V,0)$ is called a rational singularity if $p_g(V,0)=0$.

\medskip

\noindent {\bf Proposition 3.1.} [Or-Wa]\quad
Suppose $V\subseteq\c^{n+1}$ is an irreducible analytic variety,
$\sigma$ is a $\c^*$-action leaving $V$ invariant,
$$
\sigma(t,(z_0,\ldots,z_n))=(t^{q_0}z_0,\ldots, t^{q_n}z_n)
$$
and $q_i>0$ for all $i$.
Then $V$ is algebraic and the ideal of polynomials in $\c[z_0,\ldots,z_n]$
vanishing on $V$ is generated by weighted homogeneous polynomials.

Let $f(z_0,\ldots,z_n)$ be a germ of an analytic function at the origin
such that $f(0)=0$.
Suppose that $f$ has an isolated critical point at the origin.
$f$ can be developed in a convergent Taylor series
$\sum_\lambda a_\lambda z^\lambda$ where
$z^\lambda=z^{\lambda_0}_0\cdots z^{\lambda_n}_n$.
Recall that the Newton boundary $\Gamma(f)$ of $f$ is the union of compact
faces of $\Gamma_+(f)$ where $\Gamma_+(f)$ is the convex hull of the
union of the subsets $\{\lambda+(\r^+)^{n+1}\}$ for $\lambda$ such that
$a_\lambda\not= 0$.
Finally, let $\Gamma_-(f)$, the Newton polyhedron of $f$, be the cone
over $\Gamma(f)$ with vertex at $0$.
For any closed face $\Delta$ of $\Gamma(f)$, we associate the polynomial
$f_\Delta(z)=\sum_{\lambda\in\Delta}a_\lambda z^{\lambda}$.
We say that $f$ is nondegenerate if $f_\Delta$ has no critical point in
$(\c^*)^{n+1}$ for any $\Delta\in\Gamma(f)$ where $\c^*=\c-\{0\}$.
The following theorem was proved by Merle and Teissier.

\medskip

\noindent {\bf Theorem 3.2.} [Me-Te]\quad
Let $(V,0)$ be an isolated hypersurface singularity defined by a
nondegenerate holomorphic function $f:(\c^{n+1},0)\to(\c,0)$.
Then the geometric genus $p_g(V,0)=\#\{p\in\z^{n+1}\cap \Gamma_-(f):$
$p$ is positive$\}$.

Now we are ready to give the classification of 3-dimensional isolated
rational hypersurface
singularities with $\c^*$-action.

\medskip

\noindent {\bf Theorem 3.3.}\quad
Let $(V,0)$ be a 3-dimensional isolated rational hypersurface
singularity with $\c^*$-action.
Then $(V,0)$ is defined by a weighted homogeneous polynomial of one
of the 19 cases of table below such that the corresponding linear
form $\alpha$ satisfies $\alpha(1,1,1,1)>1$.

\medskip

\noindent (I) $f(x,y,z,w)=x^a+y^b+z^c+w^d$.

\smallskip
\item{1.} $(a,b,c,d)=(2,2,r,s)$, $r\geq 2$, $s\geq r$.
\smallskip
\item{2.} $(a,b,c,d)=(2,3,3,r)$, $r\geq 3$.
\smallskip
\item{3.} $(a,b,c,d)=(2,3,4,r)$, $r\geq 4$.
\smallskip
\item{4.} $(a,b,c,d)=(2,3,5,r)$, $r\geq 5$.
\smallskip
\item{5.} $(a,b,c,d)=(2,3,6,r)$, $r\geq 6$.
\smallskip
\item{6.} $(a,b,c,d)=(2,3,7,u_1)$, $u_1\in\{7,8,\ldots, 41\}$.
\smallskip
\item{7.} $(a,b,c,d)=(2,3,8,u_2)$, $u_2\in\{8,9,\ldots,23\}$.
\smallskip
\item{8.} $(a,b,c,d)=(2,3,9,u_3)$, $u_3\in\{9,10,\ldots,17\}$.
\smallskip
\item{9.} $(a,b,c,d)=(2,3,10,u_4)$, $u_4\in\{10,11,\ldots,14\}$.
\smallskip
\item{10.} $(a,b,c,d)=(2,3,11,u_5)$, $u_5\in\{11,12,13\}$.
\smallskip
\item{11.} $(a,b,c,d)=(2,4,4,r)$, $r\geq 4$.
\smallskip
\item{12.} $(a,b,c,d)=(2,4,5,u_6)$, $u_6\in\{5,6,\ldots, 19\}$.
\smallskip
\item{13.} $(a,b,c,d)=(2,4,6,u_7)$, $u_7\in\{6,7\ldots,11\}$.
\smallskip
\item{14.} $(a,b,c,d)=(2,4,7,u_8)$, $u_8\in\{7,8,9\}$.
\smallskip
\item{15.} $(a,b,c,d)=(2,5,5,u_9)$, $u_9\in\{5,6,\ldots,9\}$.
\smallskip
\item{16.} $(a,b,c,d)=(2,5,6,u_{10})$, $u_{10}\in\{6,7\}$.
\smallskip
\item{17.} $(a,b,c,d)=(3,3,3,r)$, $r\geq 3$.
\smallskip
\item{18.} $(a,b,c,d)=(3,3,4,u_{11})$, $u_{11}\in\{4,5,\ldots,11\}$.
\smallskip
\item{19.} $(a,b,c,d)=(3,3,5,u_{12})$, $u_{12}\in\{5,6,7\}$.
\smallskip
\item{20.} $(a,b,c,d)=(3,4,4,u_{13})$, $u_{13}\in\{4,5\}$.

\medskip

\noindent (II) $f(x,y,z,w)=x^a+y^b+z^c+zw^d$.

\smallskip
\item{1.}
\Itemitem{(1)} $(a,b)=(2,2)$, $\{c,d\}=\{r,s\}$, $r\geq 2$, $s\geq r$.
\itemitem{(2)} $(a,b)=(2,r)$, $\{c,d\}=\{2,s\}$, $r\geq 3$, $s\geq r$.
\itemitem{(3)} $(a,b)=(2,s)$, $\{c,d\}=\{2,r\}$, $r\geq 2$, $s\geq r+1$.
\itemitem{(4)} $(a,b)=(r,s)$, $\{c,d\}=\{2,2\}$,
where $(r,s)=(3,s)$, $s\geq 3$,
$\mathrel{\mathop{=\hskip-5pt=}\limits^{\hbox{or}}}(4,s)$,
$s\geq 4$,

\itemitem{}
$\mathrel{\mathop{=\hskip-5pt=}\limits^{\hbox{or}}}(5,u_0)$,
$u_0\in\{5,6,\ldots,19\}$,
$\mathrel{\mathop{=\hskip-5pt=}\limits^{\hbox{or}}}(6,u_1)$,
$u_1\in\{6,7,\ldots,11\}$,
$\mathrel{\mathop{=\hskip-5pt=}\limits^{\hbox{or}}}(7,u_2)$,
$u_2\in\{7,8,9\}$.

\smallskip
\item{2.}
\Itemitem{(1)} $(a,b)=(2,3)$, $\{c,d\}=\{r,s\}$,
where $(r,s)=(3,s)$, $s\geq 3$,
$\mathrel{\mathop{=\hskip-5pt=}\limits^{\hbox{or}}}(4,s)$,
$s\geq 4$,
$\mathrel{\mathop{=\hskip-5pt=}\limits^{\hbox{or}}}(5,s)$,
$s\geq 5$,
$\mathrel{\mathop{=\hskip-5pt=}\limits^{\hbox{or}}}(6,s)$,
$s\geq 6$,
$\mathrel{\mathop{=\hskip-5pt=}\limits^{\hbox{or}}}(7,u_3)$,
$u_3\in\{7,8,\ldots,35\}$,
$\mathrel{\mathop{=\hskip-5pt=}\limits^{\hbox{or}}}(8,u_4)$,
$u_4\in\{8,9,\ldots,20\}$,
$\mathrel{\mathop{=\hskip-5pt=}\limits^{\hbox{or}}}(9,u_5)$,
$u_5\in\{9,10,\ldots,15\}$,
$\mathrel{\mathop{=\hskip-5pt=}\limits^{\hbox{or}}}(10,u_6)$,
$u_6\in\{10,11,12,13\}$,
$\mathrel{\mathop{=\hskip-5pt=}\limits^{\hbox{or}}}(11,11)$.

\itemitem{(2)} $(a,b)=(2,r)$, $\{c,d\}=\{3,s\}$,
where $(r,s)=(4,s)$, $s\geq 4$,
$\mathrel{\mathop{=\hskip-5pt=}\limits^{\hbox{or}}}(5,s)$,
$s\geq 5$,
$\mathrel{\mathop{=\hskip-5pt=}\limits^{\hbox{or}}}(6,s)$,
$s\geq 6$,
$\mathrel{\mathop{=\hskip-5pt=}\limits^{\hbox{or}}}(7,u_7)$,
$u_7\in\{7,8,\ldots,27\}$,
$\mathrel{\mathop{=\hskip-5pt=}\limits^{\hbox{or}}}(8,u_8)$,
$u_8\in\{8,9,\ldots,15\}$,
$\mathrel{\mathop{=\hskip-5pt=}\limits^{\hbox{or}}}(9,u_9)$,
$u_9\in\{9,10,11\}$.

\itemitem{(3)} $(a,b)=(2,s)$, $\{c,d\}=\{3,r\}$,
where $(s,r)=(s,3)$, $s\geq 4$,
$\mathrel{\mathop{=\hskip-5pt=}\limits^{\hbox{or}}}(s,4)$,
$s\geq 5$,
$\mathrel{\mathop{=\hskip-5pt=}\limits^{\hbox{or}}}(u_{10},5)$,
$u_{10}\in\{6,7,\ldots,29\}$,
$\mathrel{\mathop{=\hskip-5pt=}\limits^{\hbox{or}}}(u_{11},6)$,
$u_{11}\in\{7,8,\ldots,17\}$,
$\mathrel{\mathop{=\hskip-5pt=}\limits^{\hbox{or}}}(u_{12},7)$,
$u_{12}\in\{8,9,\ldots,13\}$,
$\mathrel{\mathop{=\hskip-5pt=}\limits^{\hbox{or}}}(u_{13},8)$,
$u_{13}\in\{9,10,11\}$,
$\mathrel{\mathop{=\hskip-5pt=}\limits^{\hbox{or}}}(10,9)$.

\itemitem{(4)} $(a,b)=(3,r)$, $\{c,d\}=\{2,s\}$,
where $(r,s)=(3,s)$, $s\geq 3$,
$\mathrel{\mathop{=\hskip-5pt=}\limits^{\hbox{or}}}(4,s)$,
$s\geq 4$,
$\mathrel{\mathop{=\hskip-5pt=}\limits^{\hbox{or}}}(5,s)$,
$s\geq 5$,
$\mathrel{\mathop{=\hskip-5pt=}\limits^{\hbox{or}}}(6,s)$,
$s\geq 6$,
$\mathrel{\mathop{=\hskip-5pt=}\limits^{\hbox{or}}}(7,u_{14})$,
$u_{14}\in\{7,8,\ldots,20\}$,
$\mathrel{\mathop{=\hskip-5pt=}\limits^{\hbox{or}}}(8,u_{15})$,
$u_{15}\in\{8,9,10,11\}$.

\itemitem{(5)} $(a,b)=(3,s)$, $\{c,d\}=\{2,r\}$,
where $(s,r)=(s,3)$, $s\geq 4$,
$\mathrel{\mathop{=\hskip-5pt=}\limits^{\hbox{or}}}(u_{16},4)$,
$u_{16}\in\{5,6,\ldots,23\}$,
$\mathrel{\mathop{=\hskip-5pt=}\limits^{\hbox{or}}}(u_{17},5)$,
$u_{17}\in\{6,7,\ldots,14\}$,
$\mathrel{\mathop{=\hskip-5pt=}\limits^{\hbox{or}}}(u_{18},6)$,
$u_{18}\in\{7,8,\ldots,11\}$,
$\mathrel{\mathop{=\hskip-5pt=}\limits^{\hbox{or}}}(u_{19},7)$,
$u_{19}\in\{8,9,10\}$,
$\mathrel{\mathop{=\hskip-5pt=}\limits^{\hbox{or}}}(9,8)$.

\itemitem{(6)} $(a,b)=(r,s)$, $\{c,d\}=\{2,3\}$,
where $(s,r)=(4,u_{20})$, $u_{20}\in\{4,5\ldots,11\}$,
$\mathrel{\mathop{=\hskip-5pt=}\limits^{\hbox{or}}}(5,u_{21})$,
$u_{21}\in\{5,6,7\}$.

\smallskip

\item{3.}
\Itemitem{(1)} $(a,b)=(2,4)$, $\{c,d\}=\{r,s\}$,
where $(r,s)=(4,s)$, $s\geq 4$,
$\mathrel{\mathop{=\hskip-5pt=}\limits^{\hbox{or}}}(5,u_{22})$,
$u_{22}\in\{5,6,\ldots,15\}$,
$\mathrel{\mathop{=\hskip-5pt=}\limits^{\hbox{or}}}(6,u_{23})$,
$u_{23}\in\{6,7,8,9\}$,
$\mathrel{\mathop{=\hskip-5pt=}\limits^{\hbox{or}}}(7,7)$.

\itemitem{(2)} $(a,b)=(2,r)$, $\{c,d\}=\{4,s\}$,
where $(r,s)=(5,u_{24})$, $u_{24}\in\{5,6,\ldots,14\}$,
$\mathrel{\mathop{=\hskip-5pt=}\limits^{\hbox{or}}}(6,u_{25})$,
$u_{25}\in\{6,7,8\}$.

\itemitem{(3)} $(a,b)=(2,s)$, $\{c,d\}=\{4,r\}$,
where $(s,r)=(u_{26},4)$, $u_{26}\in\{5,6,\ldots,15\}$,
$\mathrel{\mathop{=\hskip-5pt=}\limits^{\hbox{or}}}(u_{27},5)$,
$u_{27}\in\{6,7,8,9\}$,
$\mathrel{\mathop{=\hskip-5pt=}\limits^{\hbox{or}}}(7,6)$.

\itemitem{(4)} $(a,b)=(4,r)$, $\{c,d\}=\{2,s\}$,
where $(r,s)=(4,s)$, $s\geq 4$,
$\mathrel{\mathop{=\hskip-5pt=}\limits^{\hbox{or}}}(5,u_{28})$,
$u_{28}\in\{5,6,\ldots,9\}$.

\itemitem{(5)} $(a,b)=(4,s)$, $\{c,d\}=\{2,r\}$,
where $(s,r)=(u_{29},4)$, $u_{29}\in\{5,6,7\}$,
$\mathrel{\mathop{=\hskip-5pt=}\limits^{\hbox{or}}}(6,5)$.

\itemitem{(6)} $(a,b)=(5,5)$, $\{c,d\}=\{2,4\}$.

\smallskip

\item{4.}
\Itemitem{(1)} $(a,b)=(2,5)$, $\{c,d\}=\{r,s\}$,
where $(r,s)=(5,u_{30})$, $u_{30}\in\{5,6,7\}$,
$\mathrel{\mathop{=\hskip-5pt=}\limits^{\hbox{or}}}(6,6)$.

\itemitem{(2)} $(a,b)=(2,u_{31})$, $u_{31}\in\{6,7\}$, $\{c,d\}=\{5,5\}$.

\smallskip

\item{5.}
\Itemitem{(1)} $(a,b)=(3,3)$, $\{c,d\}=\{r,s\}$,
where $(r,s)=(3,s)$, $s\geq 3$,
$\mathrel{\mathop{=\hskip-5pt=}\limits^{\hbox{or}}}(4,u_{32})$,
$u_{32}\in\{4,5,\ldots,8\}$,
$\mathrel{\mathop{=\hskip-5pt=}\limits^{\hbox{or}}}(5,5)$.

\itemitem{(2)} $(a,b)=(3,4)$, $\{c,d\}=\{3,u_{33}\}$,
$u_{33}\in\{4,5,6,7\}$.

\itemitem{(3)} $(a,b)=(3,s)$, $\{c,d\}=\{3,r\}$,
where $(s,r)=(u_{34},3)$, $u_{34}\in\{4,5,\ldots,8\}$,
$\mathrel{\mathop{=\hskip-5pt=}\limits^{\hbox{or}}}(5,4)$.

\itemitem{(4)} $(a,b)=(4,u_{35})$, $u_{35}\in\{4,5\}$, $\{c,d\}=\{3,3\}$.

\smallskip

\item{6.}
$(a,b)=(3,4)$, $\{c,d\}=\{4,4\}$.

\medskip

\noindent (III) $f(x,y,z,w)=x^a+y^b+z^cw+zw^d$.

\smallskip
\item{1.}
\Itemitem{(1)} $(a,b)=(2,2)$, $(c,d)=(r,s)$, $r\geq 2$, $s\geq r$.
\itemitem{(2)} $(a,b)=(2,r)$, $(c,d)=(2,s)$, $r\geq 3$, $s\geq r$.
\itemitem{(3)} $(a,b)=(2,s)$, $(c,d)=(2,r)$, $r\geq 2$, $s\geq r+1$.
\itemitem{(4)} $(a,b)=(r,s)$, $(c,d)=(2,2)$,
where $(r,s)=(3,s)$, $s\geq 3$,
$\mathrel{\mathop{=\hskip-5pt=}\limits^{\hbox{or}}}(4,u_{1})$,
$u_{1}\in\{4,5,\ldots,11\}$,
$\mathrel{\mathop{=\hskip-5pt=}\limits^{\hbox{or}}}(5,u_{2})$,
$u_{2}\in\{5,6,7\}$.

\smallskip
\item{2.}
\Itemitem{(1)} $(a,b)=(2,3)$, $(c,d)=(r,s)$,
where $(r,s)=(3,s)$, $s\geq 3$,
$\mathrel{\mathop{=\hskip-5pt=}\limits^{\hbox{or}}}(4,s)$,
$s\geq 3$,
$\mathrel{\mathop{=\hskip-5pt=}\limits^{\hbox{or}}}(5,s)$,
$s\geq 5$,
$\mathrel{\mathop{=\hskip-5pt=}\limits^{\hbox{or}}}(6,s)$,
$s\geq 6$,
$\mathrel{\mathop{=\hskip-5pt=}\limits^{\hbox{or}}}(7,u_3)$,
$u_3\in\{7,8,\ldots,30\}$,
$\mathrel{\mathop{=\hskip-5pt=}\limits^{\hbox{or}}}(8,u_4)$,
$u_4\in\{8,9,\ldots,18\}$,
$\mathrel{\mathop{=\hskip-5pt=}\limits^{\hbox{or}}}(9,u_5)$,
$u_5\in\{9,10,\ldots,14\}$,
$\mathrel{\mathop{=\hskip-5pt=}\limits^{\hbox{or}}}(10,u_6)$,
$u_6\in\{10,11,12\}$.

\itemitem{(2)} $(a,b)=(2,r)$, $(c,d)=(3,s)$,
where $(r,s)=(4,s)$, $s\geq 4$,
$\mathrel{\mathop{=\hskip-5pt=}\limits^{\hbox{or}}}(5,s)$,
$s\geq 5$,
$\mathrel{\mathop{=\hskip-5pt=}\limits^{\hbox{or}}}(6,s)$,
$s\geq 6$,
$\mathrel{\mathop{=\hskip-5pt=}\limits^{\hbox{or}}}(7,u_7)$,
$u_7\in\{7,8,\ldots,18\}$,
$\mathrel{\mathop{=\hskip-5pt=}\limits^{\hbox{or}}}(8,u_8)$,
$u_8\in\{8,9,10\}$.

\itemitem{(3)} $(a,b)=(2,s)$, $(c,d)=(3,r)$,
where $(s,r)=(s,3)$, $s\geq 4$,
$\mathrel{\mathop{=\hskip-5pt=}\limits^{\hbox{or}}}(u_9,4)$,
$u_9\in\{5,6,\ldots,21\}$,
$\mathrel{\mathop{=\hskip-5pt=}\limits^{\hbox{or}}}(u_{10},5)$,
$u_{10}\in\{6,7,\ldots,13\}$,
$\mathrel{\mathop{=\hskip-5pt=}\limits^{\hbox{or}}}(u_{11},6)$,
$u_{11}\in\{7,8,\ldots,11\}$,
$\mathrel{\mathop{=\hskip-5pt=}\limits^{\hbox{or}}}(u_{12},6)$,
$u_{12}\in\{8,9\}$,
$\mathrel{\mathop{=\hskip-5pt=}\limits^{\hbox{or}}}(9,8)$.

\itemitem{(4)} $(a,b)=(3,r)$, $(c,d)=(2,s)$,
where $(r,s)=(3,s)$, $s\geq 3$,
$\mathrel{\mathop{=\hskip-5pt=}\limits^{\hbox{or}}}(4,s)$,
$s\geq 4$,
$\mathrel{\mathop{=\hskip-5pt=}\limits^{\hbox{or}}}(5,s)$,
$s\geq 5$,
$\mathrel{\mathop{=\hskip-5pt=}\limits^{\hbox{or}}}(6,s)$,
$s\geq 6$,
$\mathrel{\mathop{=\hskip-5pt=}\limits^{\hbox{or}}}(7,u_{13})$,
$u_{13}\in\{7,8,9,10\}$.

\itemitem{(5)} $(a,b)=(3,s)$, $(c,d)=(2,r)$,
where $(s,r)=(u_{14},3)$, $u_{14}\in\{4,5,\ldots,14\}$,
$\mathrel{\mathop{=\hskip-5pt=}\limits^{\hbox{or}}}(u_{15},4)$,
$u_{15}\in\{5,6,\ldots,10\}$,
$\mathrel{\mathop{=\hskip-5pt=}\limits^{\hbox{or}}}(u_{16},5)$,
$u_{16}\in\{6,7,8\}$,
$\mathrel{\mathop{=\hskip-5pt=}\limits^{\hbox{or}}}(u_{17},6)$,
$u_{17}\in\{7,8\}$.

\itemitem{(6)} $(a,b)=(4,u_{18})$, $u_{18}\in\{4,5,6\}$, $(c,d)=(2,3)$.

\smallskip
\item{3.}
\Itemitem{(1)} $(a,b)=(2,4)$, $(c,d)=(r,s)$,
where $(r,s)=(4,s)$, $s\geq 4$,
$\mathrel{\mathop{=\hskip-5pt=}\limits^{\hbox{or}}}(5,u_{19})$,
$u_{19}\in\{5,6,\ldots,12\}$,
$\mathrel{\mathop{=\hskip-5pt=}\limits^{\hbox{or}}}(6,u_{20})$,
$u_{20}\in\{6,7,8\}$.

\itemitem{(2)} $(a,b)=(2,r)$, $(c,d)=(4,s)$,
where $(r,s)=(5,u_{21})$, $u_{21}\in\{5,6,\ldots,11\}$,
$\mathrel{\mathop{=\hskip-5pt=}\limits^{\hbox{or}}}(6,6)$.

\itemitem{(3)} $(a,b)=(2,s)$, $(c,d)=(4,r)$,
where $(s,r)=(u_{22},4)$, $u_{22}\in\{5,6,\ldots,9\}$,
$\mathrel{\mathop{=\hskip-5pt=}\limits^{\hbox{or}}}(u_{23},5)$,
$u_{23}\in\{6,7\}$.

\itemitem{(4)} $(a,b)=(4,r)$, $(c,d)=(2,s)$,
where $(r,s)=(4,s)$, $s\geq 4$,
$\mathrel{\mathop{=\hskip-5pt=}\limits^{\hbox{or}}}(5,5)$.

\itemitem{(5)} $(a,b)=(4,5)$, $(c,d)=(2,4)$.

\smallskip
\item{4.}
$(a,b)=(2,5)$, $(c,d)=(5,u_{24})$, $u_{24}\in\{5,6\}$.

\smallskip
\item{5.}
\Itemitem{(1)} $(a,b)=(3,3)$, $(c,d)=(r,s)$,
where $(r,s)=(3,s)$, $s\geq 3$,
$\mathrel{\mathop{=\hskip-5pt=}\limits^{\hbox{or}}}(4,u_{25})$,
$u_{25}\in\{4,5,6\}$.

\itemitem{(2)} $(a,b)=(3,4)$, $(c,d)=(3,u_{26})$, $u_{26}\in\{4,5\}$.

\itemitem{(3)} $(a,b)=(3,u_{27})$, $u_{27}\in\{4,5\}$, $(c,d)=(3,3)$.

\medskip

\noindent (IV) $f(x,y,z,w)=x^a+xy^b+z^c+zw^d$.

\smallskip
\item{1.}
\Itemitem{(1)} $(a,b,c,d)=(r,1,s,t)$, $r\geq 2$, $s\geq 2$, $t\geq 2$.
\itemitem{(2)} $(a,b,c,d)=(r,1,s,1)$, $r\geq 2$, $s\geq r$.

\smallskip
\item{2.}
\Itemitem{(1)} $\{a,b\}=\{2,2\}$, $\{c,d\}=\{r,s\}$,
where $(r,s)=(2,s)$, $s\geq 2$,
$\mathrel{\mathop{=\hskip-5pt=}\limits^{\hbox{or}}}(3,s)$,
$s\geq 3$,
$\mathrel{\mathop{=\hskip-5pt=}\limits^{\hbox{or}}}(4,s)$,
$s\geq 4$,
$\mathrel{\mathop{=\hskip-5pt=}\limits^{\hbox{or}}}(5,u_1)$,
$u_1\in\{5,6,\ldots,15\}$,
$\mathrel{\mathop{=\hskip-5pt=}\limits^{\hbox{or}}}(6,u_2)$,
$u_2\in\{6,7,8,9\}$,
$\mathrel{\mathop{=\hskip-5pt=}\limits^{\hbox{or}}}(7,7)$.

\itemitem{(2)} $\{a,b\}=\{2,r\}$, $\{c,d\}=\{2,s\}$, $r\geq 3$, $s\geq r$.

\smallskip
\item{3.}
\Itemitem{(1)} $\{a,b\}=\{2,3\}$, $\{c,d\}=\{r,s\}$,
where $(r,s)=(3,s)$, $s\geq 3$,
$\mathrel{\mathop{=\hskip-5pt=}\limits^{\hbox{or}}}(4,u_5)$,
$u_5\in\{4,5,\ldots,8\}$,
$\mathrel{\mathop{=\hskip-5pt=}\limits^{\hbox{or}}}(5,5)$.

\itemitem{(2)} $\{a,b\}=\{2,r\}$, $\{c,d\}=\{3,s\}$,
where $(r,s)=(4,u_6)$, $u_6\in\{4,5,\ldots,15\}$,
$\mathrel{\mathop{=\hskip-5pt=}\limits^{\hbox{or}}}(5,u_7)$,
$u_7\in\{5,6,\ldots,9\}$,
$\mathrel{\mathop{=\hskip-5pt=}\limits^{\hbox{or}}}(6,u_8)$
$u_8\in\{6,7\}$.

\itemitem{(3)} $\{a,b\}=\{2,s\}$, $\{c,d\}=\{3,r\}$,
where $(s,r)=(s,3)$, $s\geq 4$,
$\mathrel{\mathop{=\hskip-5pt=}\limits^{\hbox{or}}}(s,4)$,
$s\geq 5$,
$\mathrel{\mathop{=\hskip-5pt=}\limits^{\hbox{or}}}(u_9,5)$,
$u_9\in\{6,7,\ldots,14\}$,
$\mathrel{\mathop{=\hskip-5pt=}\limits^{\hbox{or}}}(u_{10},6)$,
$u_{10}\in\{7,8\}$.

\smallskip
\item{4.}
\Itemitem{(1)} $\{a,b\}=\{2,4\}$, $\{c,d\}=\{4,u_{16}\}$, $u_{16}\in\{4,5\}$.
\itemitem{(2)} $\{a,b\}=\{2,u_{17}\}$, $u_{17}\in\{5,6,7\}$, $\{c,d\}=\{4,4\}$.

\smallskip
\item{5.}
$\{a,b\}=\{3,3\}$, $\{c,d\}=\{3,u_{19}\}$, $u_{19}\in\{3,4,5\}$.

\medskip

\noindent (V) $f(x,y,z,w)=x^ay+xy^b+z^c+zw^d$.

\smallskip
\item{1.}
$(a,b,c,d)=(r,s,t,1)$, $r\geq 2$, $s\geq r$, $t\geq 2$.

\smallskip
\item{2.}
\Itemitem{(1)} $(a,b)=(2,2)$, $\{c,d\}=\{r,s\}$,
where $(r,s)=(2,s)$, $s\geq 2$,
$\mathrel{\mathop{=\hskip-5pt=}\limits^{\hbox{or}}}(3,s)$,
$s\geq 3$,
$\mathrel{\mathop{=\hskip-5pt=}\limits^{\hbox{or}}}(4,u_1)$,
$u_1\in\{4,5,\ldots,8\}$,
$\mathrel{\mathop{=\hskip-5pt=}\limits^{\hbox{or}}}(5,5)$.

\itemitem{(2)} $(a,b)=(2,r)$, $\{c,d\}=\{2,s\}$, $r\geq 3$, $s\geq r$.

\itemitem{(3)} $(a,b)=(2,s)$, $\{c,d\}=\{2,r\}$, $r\geq 2$, $s\geq r+1$.
\itemitem{(4)} $(a,b)=(r,s)$, $\{c,d\}=\{2,2\}$,
where $(r,s)=(3,s)$, $s\geq 3$,
$\mathrel{\mathop{=\hskip-5pt=}\limits^{\hbox{or}}}(4,s)$,
$s\geq 4$,
$\mathrel{\mathop{=\hskip-5pt=}\limits^{\hbox{or}}}(5,u_2)$,
$u_2\in\{5,6,\ldots,12\}$,
$\mathrel{\mathop{=\hskip-5pt=}\limits^{\hbox{or}}}(6,u_3)$,
$u_3\in\{6,7,8\}$.

\smallskip

\item{3.}
\Itemitem{(1)} $(a,b)=(2,3)$, $\{c,d\}=\{r,s\}$,
where $(r,s)=(3,u_4)$, $u_4\in\{3,4\ldots,9\}$,
$\mathrel{\mathop{=\hskip-5pt=}\limits^{\hbox{or}}}(4,4)$.

\itemitem{(2)} $(a,b)=(2,r)$, $\{c,d\}=\{3,s\}$,
where $(r,s)=(4,u_5)$, $u_5\in\{4,5,6\}$,
$\mathrel{\mathop{=\hskip-5pt=}\limits^{\hbox{or}}}(5,5)$.

\itemitem{(3)} $(a,b)=(2,s)$, $\{c,d\}=\{3,r\}$,
where $(s,r)=(s,3)$, $s\geq 4$,
$\mathrel{\mathop{=\hskip-5pt=}\limits^{\hbox{or}}}(s,4)$,
$s\geq 5$,
$\mathrel{\mathop{=\hskip-5pt=}\limits^{\hbox{or}}}(u_6,5)$,
$u_6\in\{6,7\}$.

\itemitem{(4)} $(a,b)=(3,r)$, $\{c,d\}=\{2,s\}$,
where $(s,r)=(3,s)$, $s\geq 3$,
$\mathrel{\mathop{=\hskip-5pt=}\limits^{\hbox{or}}}(4,u_7)$,
$u_7\in\{4,5,\ldots,10\}$,
$\mathrel{\mathop{=\hskip-5pt=}\limits^{\hbox{or}}}(5,u_8)$,
$u_8\in\{5,6\}$.

\itemitem{(5)} $(a,b)=(3,s)$, $\{c,d\}=\{2,r\}$,
where $(s,r)=(s,3)$, $s\geq 4$,
$\mathrel{\mathop{=\hskip-5pt=}\limits^{\hbox{or}}}(u_9,4)$,
$u_9\in\{5,6,\ldots,10\}$,
$\mathrel{\mathop{=\hskip-5pt=}\limits^{\hbox{or}}}(6,5)$.

\itemitem{(6)} $(a,b)=(4,u_{10})$, $u_{10}\in\{4,5,6\}$, $\{c,d\}=\{2,3\}$.

\smallskip

\item{4.}
\Itemitem{(1)} $(a,b)=(2,4)$, $\{c,d\}=\{4,4\}$.
\itemitem{(2)} $(a,b)=(4,4)$, $\{c,d\}=\{2,4\}$.

\smallskip

\item{5.}
\Itemitem{(1)} $(a,b)=(3,3)$, $\{c,d\}=\{3,3\}$.
\itemitem{(2)} $(a,b)=(3,4)$, $\{c,d\}=\{3,3\}$.

\medskip

\noindent (VI) $f(x,y,z,w)=x^ay+xy^b+z^cw+zw^d$.

\smallskip
\item{1.}
\Itemitem{(1)} $(a,b)=(2,2)$, $(c,d)=(r,s)$,
where $(r,s)=(2,s)$, $s\geq 2$,
$\mathrel{\mathop{=\hskip-5pt=}\limits^{\hbox{or}}}(3,s)$,
$s\geq 3$,
$\mathrel{\mathop{=\hskip-5pt=}\limits^{\hbox{or}}}(4,u_1)$,
$u_1\in\{4,5,6\}$.

\itemitem{(2)} $(a,b)=(2,r)$, $(c,d)=(2,s)$, $r\geq 3$, $s\geq r$.

\smallskip
\item{2.}
\Itemitem{(1)} $(a,b)=(2,3)$, $(c,d)=(3,u_3)$, $u_3\in\{3,4,5,6\}$.
\itemitem{(2)} $(a,b)=(2,4)$, $(c,d)=(3,4)$.
\itemitem{(3)} $(a,b)=(2,s)$, $(c,d)=(3,r)$,
where $(s,r)=(s,3)$, $s\geq 4$,
$\mathrel{\mathop{=\hskip-5pt=}\limits^{\hbox{or}}}(5,4)$.

\medskip

\noindent (VII) $f(x,y,z,w)=x^a+y^b+yz^c+zw^d$.

\smallskip
\item{1.}
\Itemitem{(1)} $(a,b,c,d)=(r,s,1,t)$, $r\geq 2$, $s\geq 2$, $t\geq 2$.
\itemitem{(2)} $(a,b,c,d)=(r,s,t,1)$, $r\geq 2$, $s\geq 2$, $t\geq 2$.
\itemitem{(3)} $(a,b,c,d)=(r,s,1,1)$, $r\geq 2$, $s\geq 2$.

\smallskip
\item{2.}
\Itemitem{(1)} $\{a,c\}=\{2,2\}$, $\{b,d\}=\{r,s\}$, $r\geq 2$, $s\geq 2$.
\itemitem{(2)} $(a,c)=(2,r)$, $\{b,d\}=\{2,s\}$, $r\geq 3$, $s\geq r$.
\itemitem{(3)} $(a,c)=(r,2)$, $\{b,d\}=\{2,s\}$,
where $(r,s)=(3,s)$, $s\geq 3$,
$\mathrel{\mathop{=\hskip-5pt=}\limits^{\hbox{or}}}(4,s)$,
$s\geq 4$,
$\mathrel{\mathop{=\hskip-5pt=}\limits^{\hbox{or}}}(5,u_1)$,
$u_1\in\{5,6,\ldots,14\}$,
$\mathrel{\mathop{=\hskip-5pt=}\limits^{\hbox{or}}}(6,u_2)$,
$u_2\in\{6,7,8\}$.
\itemitem{(4)} $(a,c)=(2,s)$, $\{b,d\}=\{2,r\}$, $r\geq 2$, $s\geq r+1$.
\itemitem{(5)} $(a,c)=(s,2)$, $\{b,d\}=\{2,r\}$,
where $(s,r)=(s,2)$, $s\geq 3$,
$\mathrel{\mathop{=\hskip-5pt=}\limits^{\hbox{or}}}(s,3)$,
$s\geq 4$,
$\mathrel{\mathop{=\hskip-5pt=}\limits^{\hbox{or}}}(u_3,4)$,
$u_3\in\{5,6,\ldots,15\}$,
$\mathrel{\mathop{=\hskip-5pt=}\limits^{\hbox{or}}}(u_4,5)$,
$u_4\in\{6,7,8,9\}$,
$\mathrel{\mathop{=\hskip-5pt=}\limits^{\hbox{or}}}(7,6)$.
\itemitem{(6)} $(a,c)=\{r,s\}$, $\{b,d\}=\{2,2\}$, $r\geq 3$, $s\geq r$.

\smallskip
\item{3.}
\Itemitem{(1)} $(a,c)=(2,3)$, $\{b,d\}=\{r,s\}$,
where $(r,s)=(3,s)$, $s\geq 3$,
$\mathrel{\mathop{=\hskip-5pt=}\limits^{\hbox{or}}}(4,s)$,
$s\geq 4$,
$\mathrel{\mathop{=\hskip-5pt=}\limits^{\hbox{or}}}(5,u_5)$,
$u_5\in\{5,6,\ldots,21\}$,
$\mathrel{\mathop{=\hskip-5pt=}\limits^{\hbox{or}}}(6,u_6)$,
$u_6\in\{6,7,\ldots,12\}$,
$\mathrel{\mathop{=\hskip-5pt=}\limits^{\hbox{or}}}(7,u_7)$,
$u_7\in\{7,8,9\}$,
$\mathrel{\mathop{=\hskip-5pt=}\limits^{\hbox{or}}}(8,8)$.

\itemitem{(2)} $(a,c)=(3,2)$, $\{b,d\}=\{r,s\}$,
where $(r,s)=(3,s)$, $s\geq 3$,
$\mathrel{\mathop{=\hskip-5pt=}\limits^{\hbox{or}}}(4,u_8)$,
$u_8\in\{4,5,\ldots,14\}$,
$\mathrel{\mathop{=\hskip-5pt=}\limits^{\hbox{or}}}(5,u_9)$,
$u_9\in\{5,6,7,8\}$,
$\mathrel{\mathop{=\hskip-5pt=}\limits^{\hbox{or}}}(6,6)$.

\itemitem{(3)} $(a,c)=(2,r)$, $\{b,d\}=\{3,s\}$,
where $(r,s)=(4,s)$, $s\geq 4$,
$\mathrel{\mathop{=\hskip-5pt=}\limits^{\hbox{or}}}(5,u_{10})$,
$u_{10}\in\{5,6,\ldots,25\}$,
$\mathrel{\mathop{=\hskip-5pt=}\limits^{\hbox{or}}}(6,u_{11})$,
$u_{11}\in\{6,7,\ldots,15\}$,
$\mathrel{\mathop{=\hskip-5pt=}\limits^{\hbox{or}}}(7,u_{12})$,
$u_{12}\in\{7,8,\ldots,12\}$,
$\mathrel{\mathop{=\hskip-5pt=}\limits^{\hbox{or}}}(8,u_{13})$,
$u_{13}\in\{8,9,10\}$,
$\mathrel{\mathop{=\hskip-5pt=}\limits^{\hbox{or}}}(9,9)$.

\itemitem{(4)} $(a,c)=(r,2)$, $\{b,d\}=\{3,s\}$,
where $(r,s)=(4,u_{14})$, $u_{14}\in\{4,5,6,7\}$.

\itemitem{(5)} $(a,c)=(2,s)$, $\{b,d\}=\{3,r\}$,
where $(s,r)=(s,3)$, $s\geq 4$,
$\mathrel{\mathop{=\hskip-5pt=}\limits^{\hbox{or}}}(s,4)$,
$s\geq 5$,
$\mathrel{\mathop{=\hskip-5pt=}\limits^{\hbox{or}}}(s,5)$,
$s\geq 6$,
$\mathrel{\mathop{=\hskip-5pt=}\limits^{\hbox{or}}}(s,6)$,
$s\geq 7$,
$\mathrel{\mathop{=\hskip-5pt=}\limits^{\hbox{or}}}(u_{14},7)$,
$u_{14}\in\{8,9,\ldots,23\}$,
$\mathrel{\mathop{=\hskip-5pt=}\limits^{\hbox{or}}}(u_{15},8)$,
$u_{15}\in\{9,10,\ldots,13\}$,
$\mathrel{\mathop{=\hskip-5pt=}\limits^{\hbox{or}}}(10,9)$.

\itemitem{(6)} $(a,c)=(s,2)$, $\{b,d\}=\{3,r\}$,
where $(s,r)=(u_{16},3)$, $u_{16}\in\{4,5,\ldots,8\}$,
$\mathrel{\mathop{=\hskip-5pt=}\limits^{\hbox{or}}}(5,4)$.

\itemitem{(7)} $(a,c)=(3,r)$, $\{b,d\}=\{2,s\}$,
where $(r,s)=(3,s)$, $s\geq 3$,
$\mathrel{\mathop{=\hskip-5pt=}\limits^{\hbox{or}}}(4,u_{17})$,
$u_{17}\in\{4,5,\ldots,20\}$,
$\mathrel{\mathop{=\hskip-5pt=}\limits^{\hbox{or}}}(5,u_{18})$,
$u_{18}\in\{5,6,\ldots,13\}$,
$\mathrel{\mathop{=\hskip-5pt=}\limits^{\hbox{or}}}(6,u_{19})$,
$u_{19}\in\{6,7,\ldots,10\}$,
$\mathrel{\mathop{=\hskip-5pt=}\limits^{\hbox{or}}}(7,u_{20})$,
$u_{20}\in\{7,8,9\}$,
$\mathrel{\mathop{=\hskip-5pt=}\limits^{\hbox{or}}}(8,8)$.

\itemitem{(8)} $(a,c)=(r,3)$, $\{b,d\}=\{2,s\}$,
where $(r,s)=(4,u_{21})$, $u_{21}\in\{4,5,\ldots,9\}$,
$\mathrel{\mathop{=\hskip-5pt=}\limits^{\hbox{or}}}(5,u_{22})$,
$u_{22}\in\{5,6\}$.

\itemitem{(9)} $(a,c)=(3,s)$, $\{b,d\}=\{2,r\}$,
where $(s,r)=(s,3)$, $s\geq 4$,
$\mathrel{\mathop{=\hskip-5pt=}\limits^{\hbox{or}}}(s,4)$,
$s\geq 5$,
$\mathrel{\mathop{=\hskip-5pt=}\limits^{\hbox{or}}}(s,5)$,
$s\geq 6$,
$\mathrel{\mathop{=\hskip-5pt=}\limits^{\hbox{or}}}(s,6)$,
$s\geq 7$,
$\mathrel{\mathop{=\hskip-5pt=}\limits^{\hbox{or}}}(u_{23},7)$,
$u_{23}\in\{8,9,\ldots,17\}$,
$\mathrel{\mathop{=\hskip-5pt=}\limits^{\hbox{or}}}(u_{24},8)$,
$u_{24}\in\{9,10\}$.

\itemitem{(10)} $(a,c)=(s,3)$, $\{b,d\}=\{2,r\}$,
where $(s,r)=(u_{25},3)$, $u_{25}\in\{4,5,\ldots,17\}$,
$\mathrel{\mathop{=\hskip-5pt=}\limits^{\hbox{or}}}(u_{26},4)$,
$u_{26}\in\{5,6,7\}$.

\itemitem{(11)} $(a,c)=(r,s)$, $\{b,d\}=\{2,3\}$,
where $(r,s)=(4,s)$, $s\geq 4$,
$\mathrel{\mathop{=\hskip-5pt=}\limits^{\hbox{or}}}(5,s)$,
$s\geq 5$,
$\mathrel{\mathop{=\hskip-5pt=}\limits^{\hbox{or}}}(6,s)$,
$s\geq 6$,
$\mathrel{\mathop{=\hskip-5pt=}\limits^{\hbox{or}}}(7,u_{27})$,
$u_{27}\in\{7,8,\ldots,13\}$.

\itemitem{(12)} $(a,c)=(s,r)$, $\{b,d\}=\{2,3\}$,
where $(s,r)=(u_{28},4)$, $u_{28}\in\{5,6,\ldots,11\}$,
$\mathrel{\mathop{=\hskip-5pt=}\limits^{\hbox{or}}}(u_{29},5)$,
$u_{29}\in\{6,7,8,9\}$,
$\mathrel{\mathop{=\hskip-5pt=}\limits^{\hbox{or}}}(u_{30},6)$,
$u_{30}\in\{7,8\}$,
$\mathrel{\mathop{=\hskip-5pt=}\limits^{\hbox{or}}}(8,7)$.

\smallskip
\item{4.}
\Itemitem{(1)} $(a,c)=(2,4)$, $\{b,d\}=\{r,s\}$,
where $(r,s)=(4,u_{31})$, $u_{31}\in\{4,5,\ldots,12\}$,
$\mathrel{\mathop{=\hskip-5pt=}\limits^{\hbox{or}}}(5,u_{32})$,
$u_{32}\in\{5,6,7\}$,
$\mathrel{\mathop{=\hskip-5pt=}\limits^{\hbox{or}}}(6,6)$.

\itemitem{(2)} $(a,c)=(4,2)$, $\{b,d\}=\{4,4\}$.

\itemitem{(3)} $(a,c)=(2,r)$, $\{b,d\}=\{4,s\}$,
where $(r,s)=(5,u_{33})$, $u_{33}\in\{5,6,7,8\}$,
$\mathrel{\mathop{=\hskip-5pt=}\limits^{\hbox{or}}}(6,6)$.

\itemitem{(4)} $(a,c)=(2,s)$, $\{b,d\}=\{4,r\}$,
where $(s,r)=(s,4)$, $s\geq 5$,
$\mathrel{\mathop{=\hskip-5pt=}\limits^{\hbox{or}}}(u_{34},5)$,
$u_{34}\in\{6,7,\ldots,11\}$,
$\mathrel{\mathop{=\hskip-5pt=}\limits^{\hbox{or}}}(7,6)$.

\itemitem{(5)} $(a,c)=(4,r)$, $\{b,d\}=\{2,s\}$,
where $(r,s)=(4,u_{35})$, $u_{35}\in\{4,5,6\}$,
$\mathrel{\mathop{=\hskip-5pt=}\limits^{\hbox{or}}}(5,5)$.

\itemitem{(6)} $(a,c)=(4,s)$, $\{b,d\}=\{2,r\}$,
where $(s,r)=(s,4)$, $s\geq 5$,
$\mathrel{\mathop{=\hskip-5pt=}\limits^{\hbox{or}}}(u_{36},5)$,
$u_{36}\in\{6,7\}$.

\itemitem{(7)} $(a,c)=(s,4)$, $\{b,d\}=\{2,r\}$,
where $(s,r)=(u_{37},4)$, $u_{37}\in\{5,6\}$.

\itemitem{(8)} $(a,c)=(r,s)$, $\{b,d\}=\{2,4\}$,
where $(r,s)=(5,u_{38})$, $u_{38}\in\{5,6,7\}$.

\smallskip
\item{5.}
\Itemitem{(1)} $(a,c)=(2,5)$, $\{b,d\}=\{5,5\}$.
\itemitem{(2)} $(a,c)=(2,6)$, $\{b,d\}=\{5,5\}$.

\smallskip
\item{6.}
\Itemitem{(1)} $(a,c)=(3,3)$, $\{b,d\}=\{r,s\}$,
where $(r,s)=(3,u_{39})$, $u_{39}\in\{3,4,5,6\}$,
$\mathrel{\mathop{=\hskip-5pt=}\limits^{\hbox{or}}}(4,4)$.
\itemitem{(2)} $(a,c)=(3,4)$, $\{b,d\}=\{3,4\}$.

\itemitem{(3)} $(a,c)=(3,s)$, $\{b,d\}=\{3,r\}$,
where $(s,r)=(s,3)$, $s\geq 4$,
$\mathrel{\mathop{=\hskip-5pt=}\limits^{\hbox{or}}}(5,4)$.

\itemitem{(4)} $(a,c)=(s,3)$, $\{b,d\}=\{3,r\}$,
where $(s,r)=(u_{40},3)$, $u_{40}\in\{4,5\}$.

\itemitem{(5)} $(a,c)=(r,s)$, $\{b,d\}=\{3,3\}$,
where $(r,s)=(4,u_{41})$, $u_{41}\in\{4,5\}$.

\medskip

\noindent (VIII) $f(x,y,z,w)=x^a+y^b+yz^c+yw^d+z^pw^q$,
$\frac{p(b-1)}{bc}+\frac{q(b-1)}{bd}=1$.

\smallskip
\item{1.}
\Itemitem{(1)} $(a,b,c,d)=(r,s,1,t)$, $r\geq 2$, $s\geq 2$, $t\geq 2$.

\itemitem{(2)} $(a,b,c,d)=(r,s,t,1)$, $r\geq 2$, $s\geq 2$, $t\geq 2$.

\itemitem{(3)} $(a,b,c,d)=(r,s,1,1)$, $r\geq 2$, $s\geq 2$.

\smallskip
\item{2.}
\Itemitem{(1)} $\{a,b\}=\{2,2\}$, $\{c,d\}=\{r,s\}$, $r\geq 2$, $s\geq r$.
\itemitem{(2)} $(a,b)=(2,r)$, $\{c,d\}=\{2,s\}$, $r\geq 3$, $s\geq r$.
\itemitem{(3)} $(a,b)=(r,2)$, $\{c,d\}=\{2,s\}$,
where $(r,s)=(3,s)$, $s\geq 3$,
$\mathrel{\mathop{=\hskip-5pt=}\limits^{\hbox{or}}}(4,s)$,
$s\geq 4$,
$\mathrel{\mathop{=\hskip-5pt=}\limits^{\hbox{or}}}(5,u_1)$,
$u_1\in\{5,6,\ldots,9\}$.

\itemitem{(4)} $(a,b)=(2,s)$, $\{c,d\}=\{2,r\}$, $r\geq 2$, $s\geq r+1$.

\itemitem{(5)} $(a,b)=(s,2)$, $\{c,d\}=\{2,r\}$,
where $(s,r)=(s,2)$, $s\geq 3$,
$\mathrel{\mathop{=\hskip-5pt=}\limits^{\hbox{or}}}(u_2,3)$,
$u_2\in\{4,5,\ldots,11\}$,
$\mathrel{\mathop{=\hskip-5pt=}\limits^{\hbox{or}}}(u_3,4)$,
$u_3\in\{5,6,7\}$,
$\mathrel{\mathop{=\hskip-5pt=}\limits^{\hbox{or}}}(6,5)$.

\itemitem{(6)} $\{a,b\}=\{r,s\}$, $\{c,d\}=\{2,2\}$, $r\geq 3$, $s\geq r$.

\smallskip
\item{3.}
\Itemitem{(1)} $(a,b)=(2,3)$, $\{c,d\}=\{r,s\}$,
where $(r,s)=(3,s)$, $s\geq 3$,
$\mathrel{\mathop{=\hskip-5pt=}\limits^{\hbox{or}}}(4,s)$,
$s\geq 4$,
$\mathrel{\mathop{=\hskip-5pt=}\limits^{\hbox{or}}}(5,u_4)$,
$u_4\in\{5,6,\ldots,19\}$,
$\mathrel{\mathop{=\hskip-5pt=}\limits^{\hbox{or}}}(6,u_5)$,
$u_5\in\{6,7,\ldots,11\}$,
$\mathrel{\mathop{=\hskip-5pt=}\limits^{\hbox{or}}}(7,u_6)$,
$u_6\in\{7,8,9\}$.

\itemitem{(2)} $(a,b)=(3,2)$, $\{c,d\}=\{r,s\}$,
where $(r,s)=(3,s)$, $s\geq 3$,
$\mathrel{\mathop{=\hskip-5pt=}\limits^{\hbox{or}}}(4,u_7)$,
$u_7\in\{4,5,\ldots,11\}$,
$\mathrel{\mathop{=\hskip-5pt=}\limits^{\hbox{or}}}(5,u_8)$,
$u_8\in\{5,6,7\}$.

\itemitem{(3)} $(a,b)=(2,r)$, $\{c,d\}=\{3,s\}$,
where $(r,s)=(4,s)$, $s\geq 4$,
$\mathrel{\mathop{=\hskip-5pt=}\limits^{\hbox{or}}}(5,u_9)$,
$u_9\in\{5,6,\ldots,23\}$,
$\mathrel{\mathop{=\hskip-5pt=}\limits^{\hbox{or}}}(6,u_{10})$,
$u_{10}\in\{6,7,\ldots,14\}$,
$\mathrel{\mathop{=\hskip-5pt=}\limits^{\hbox{or}}}(7,u_{11})$,
$u_{11}\in\{7,8,\ldots,11\}$,
$\mathrel{\mathop{=\hskip-5pt=}\limits^{\hbox{or}}}(8,u_{12})$,
$u_{12}\in\{8,9,10\}$,
$\mathrel{\mathop{=\hskip-5pt=}\limits^{\hbox{or}}}(9,9)$.

\itemitem{(4)} $(a,b)=(r,2)$, $\{c,d\}=\{3,s\}$,
where $(r,s)=(4,u_{13})$, $u_{13}\in\{4,5\}$.

\itemitem{(5)} $(a,b)=(2,s)$, $\{c,d\}=\{3,r\}$,
where $(s,r)=(s,3)$, $s\geq 4$,
$\mathrel{\mathop{=\hskip-5pt=}\limits^{\hbox{or}}}(s,4)$,
$s\geq 5$,
$\mathrel{\mathop{=\hskip-5pt=}\limits^{\hbox{or}}}(s,5)$,
$s\geq 6$,
$\mathrel{\mathop{=\hskip-5pt=}\limits^{\hbox{or}}}(s,6)$,
$s\geq 7$,
$\mathrel{\mathop{=\hskip-5pt=}\limits^{\hbox{or}}}(u_{14},7)$,
$u_{14}\in\{8,9,\ldots,21\}$,
$\mathrel{\mathop{=\hskip-5pt=}\limits^{\hbox{or}}}(u_{15},8)$,
$u_{15}\in\{9,10,11,12\}$.

\itemitem{(6)} $(a,b)=(s,2)$, $\{c,d\}=\{3,r\}$,
where $(s,r)=(u_{16},3)$, $u_{16}\in\{4,5\}$.

\itemitem{(7)} $(a,b)=(3,r)$, $\{c,d\}=\{2,s\}$,
where $(r,s)=(3,s)$, $s\geq 3$,
$\mathrel{\mathop{=\hskip-5pt=}\limits^{\hbox{or}}}(4,u_{17})$,
$u_{17}\in\{4,5,\ldots,17\}$,
$\mathrel{\mathop{=\hskip-5pt=}\limits^{\hbox{or}}}(5,u_{18})$,
$u_{18}\in\{5,6,\ldots,11\}$,
$\mathrel{\mathop{=\hskip-5pt=}\limits^{\hbox{or}}}(6,u_{19})$,
$u_{19}\in\{6,7,8,9\}$,
$\mathrel{\mathop{=\hskip-5pt=}\limits^{\hbox{or}}}(7,u_{20})$,
$u_{20}\in\{7,8\}$,
$\mathrel{\mathop{=\hskip-5pt=}\limits^{\hbox{or}}}(8,8)$.

\itemitem{(8)} $(a,b)=(r,3)$, $\{c,d\}=\{2,s\}$,
where $(r,s)=(4,u_{21})$, $u_{21}\in\{4,5,6,7\}$.

\itemitem{(9)} $(a,b)=(3,s)$, $\{c,d\}=\{2,r\}$,
where $(s,r)=(s,3)$, $s\geq 4$,
$\mathrel{\mathop{=\hskip-5pt=}\limits^{\hbox{or}}}(s,4)$,
$s\geq 5$,
$\mathrel{\mathop{=\hskip-5pt=}\limits^{\hbox{or}}}(s,5)$,
$s\geq 6$,
$\mathrel{\mathop{=\hskip-5pt=}\limits^{\hbox{or}}}(s,6)$,
$s\geq 7$,
$\mathrel{\mathop{=\hskip-5pt=}\limits^{\hbox{or}}}(u_{22},7)$,
$u_{22}\in\{8,9,\ldots,14\}$.

\itemitem{(10)} $(a,b)=(s,3)$, $\{c,d\}=\{2,r\}$,
where $(s,r)=(u_{23},3)$, $u_{23}\in\{4,5,\ldots,8\}$,
$\mathrel{\mathop{=\hskip-5pt=}\limits^{\hbox{or}}}(5,4)$.

\itemitem{(11)} $(a,b)=(r,s)$, $\{c,d\}=\{2,3\}$,
where $(r,s)=(4,s)$, $s\geq 4$,
$\mathrel{\mathop{=\hskip-5pt=}\limits^{\hbox{or}}}(5,s)$,
$s\geq 5$,
$\mathrel{\mathop{=\hskip-5pt=}\limits^{\hbox{or}}}(6,s)$,
$s\geq 6$.

\itemitem{(12)} $(a,b)=(s,r)$, $\{c,d\}=\{2,3\}$,
where $(s,r)=(u_{24},4)$, $u_{24}\in\{5,6,7\}$,
$\mathrel{\mathop{=\hskip-5pt=}\limits^{\hbox{or}}}(u_{25},5)$,
$u_{25}\in\{6,7\}$,
$\mathrel{\mathop{=\hskip-5pt=}\limits^{\hbox{or}}}(7,6)$.

\smallskip
\item{4.}
\Itemitem{(1)} $(a,b)=(2,4)$, $\{c,d\}=\{r,s\}$,
where $(r,s)=(4,u_{26})$, $u_{26}\in\{4,5,\ldots,11\}$,
$\mathrel{\mathop{=\hskip-5pt=}\limits^{\hbox{or}}}(5,u_{27})$,
$u_{27}\in\{5,6,7\}$.

\itemitem{(2)} $(a,b)=(2,r)$, $\{c,d\}=\{4,s\}$,
where $(r,s)=(5,u_{28})$, $u_{28}\in\{5,6,7\}$,
$\mathrel{\mathop{=\hskip-5pt=}\limits^{\hbox{or}}}(6,6)$.

\itemitem{(3)} $(a,b)=(2,s)$, $\{c,d\}=\{4,r\}$,
where $(s,r)=(s,4)$, $s\geq 5$,
$\mathrel{\mathop{=\hskip-5pt=}\limits^{\hbox{or}}}(u_{29},5)$
$u_{29}\in\{6,7,\ldots,10\}$.

\itemitem{(4)} $(a,b)=(4,r)$, $\{c,d\}=\{2,s\}$,
where $(r,s)=(4,u_{30})$, $u_{30}\in\{4,5\}$,
$\mathrel{\mathop{=\hskip-5pt=}\limits^{\hbox{or}}}(5,5)$.

\itemitem{(5)} $(a,b)=(4,s)$, $\{c,d\}=\{2,r\}$,
where $(s,r)=(s,4)$, $s\geq 5$.

\itemitem{(6)} $(a,b)=(5,4)$, $\{c,d\}=\{2,4\}$.

\smallskip
\item{5.}
$(a,b)=(2,5)$, $\{c,d\}=\{5,5\}$.

\smallskip
\item{6.}
\Itemitem{(1)} $(a,b)=(3,3)$, $\{c,d\}=\{r,s\}$,
where $(r,s)=(3,u_{31})$, $u_{31}\in\{3,4,5\}$.

\itemitem{(2)} $(a,b)=(3,4)$, $\{c,d\}=\{3,4\}$.

\itemitem{(3)} $(a,b)=(3,s)$, $\{c,d\}=\{3,r\}$,
where $(s,r)=(s,3)$, $s\geq 4$.

\itemitem{(4)} $(a,b)=(4,3)$, $\{c,d\}=\{3,3\}$.

\medskip

\noindent (IX) $f(x,y,z,w)=x^a+y^bw+z^cw+yw^d+y^pz^q$,
$\frac{p(d-1)}{bd-1}+\frac{qb(d-1)}{c(bd-1)}=1$.

\smallskip
\item{1.}
$(a,b,c,d)=(r,s,1,t)$, $r\geq 2$, $s\geq 2$, $t\geq 2$.

\smallskip
\item{2.}
\Itemitem{(1)} $\{a,b\}=\{2,2\}$, $\{c,d\}=\{r,s\}$, $r\geq 2$, $s\geq r$.
\itemitem{(2)} $(a,b)=(2,r)$, $\{c,d\}=\{2,s\}$, $r\geq 3$, $s\geq r$.
\itemitem{(3)} $(a,b,c,d)=(r,2,2,s)$, $r\geq 3$, $s\geq r$.

\itemitem{(4)} $(a,b,c,d)=(r,2,s,2)$,
where $(r,s)=(3,s)$, $s\geq 3$,
$\mathrel{\mathop{=\hskip-5pt=}\limits^{\hbox{or}}}(4,u_1)$,
$u_1\in\{4,5,6,7\}$.

\itemitem{(5)} $(a,b)=(2,s)$, $\{c,d\}=\{2,r\}$, $r\geq 2$, $s\geq r+1$.

\itemitem{(6)} $(a,b,c,d)=(s,2,2,r)$, $r\geq 2$, $s\geq r+1$.

\itemitem{(7)} $(a,b,c,d)=(s,2,r,2)$,
where $(s,r)=(u_2,3)$, $u_2\in\{4,5,\ldots,8\}$,
$\mathrel{\mathop{=\hskip-5pt=}\limits^{\hbox{or}}}(5,4)$.

\itemitem{(8)} $(a,b,c,d)=(r,s,2,2)$,
where $(r,s)=(3,s)$, $s\geq 3$,
$\mathrel{\mathop{=\hskip-5pt=}\limits^{\hbox{or}}}(4,s)$,
$s\geq 4$,
$\mathrel{\mathop{=\hskip-5pt=}\limits^{\hbox{or}}}(5,u_3)$,
$u_3\in\{5,6,7\}$.

\itemitem{(9)} $(a,b,c,d)=(s,r,2,2)$,
where $(s,r)=(u_4,3)$, $u_4\in\{4,5,\ldots,9\}$,
$\mathrel{\mathop{=\hskip-5pt=}\limits^{\hbox{or}}}(u_5,4)$,
$u_5\in\{5,6\}$.

\smallskip

\item{3.}
\Itemitem{(1)} $(a,b,c,d)=(2,3,r,s)$,
where $(r,s)=(3,s)$, $s\geq 3$,
$\mathrel{\mathop{=\hskip-5pt=}\limits^{\hbox{or}}}(4,s)$,
$s\geq 4$,
$\mathrel{\mathop{=\hskip-5pt=}\limits^{\hbox{or}}}(5,s)$,
$s\geq 5$,
$\mathrel{\mathop{=\hskip-5pt=}\limits^{\hbox{or}}}(6,s)$,
$s\geq 6$,
$\mathrel{\mathop{=\hskip-5pt=}\limits^{\hbox{or}}}(7,u_7)$,
$u_7\in\{7,8,\ldots,14\}$,
$\mathrel{\mathop{=\hskip-5pt=}\limits^{\hbox{or}}}(8,8)$.

\itemitem{(2)} $(a,b,c,d)=(2,3,s,r)$,
where $(s,r)=(s,3)$, $s\geq 4$,
$\mathrel{\mathop{=\hskip-5pt=}\limits^{\hbox{or}}}(u_8,4)$,
$u_8\in\{5,6,\ldots,17\}$,
$\mathrel{\mathop{=\hskip-5pt=}\limits^{\hbox{or}}}(u_9,5)$,
$u_9\in\{6,7,\ldots,11\}$,
$\mathrel{\mathop{=\hskip-5pt=}\limits^{\hbox{or}}}(u_{10},6)$,
$u_{10}\in\{7,8,9\}$,
$\mathrel{\mathop{=\hskip-5pt=}\limits^{\hbox{or}}}(8,7)$.

\itemitem{(3)} $(a,b,c,d)=(3,2,r,s)$,
where $(r,s)=(3,s)$, $s\geq 3$,
$\mathrel{\mathop{=\hskip-5pt=}\limits^{\hbox{or}}}(4,s)$,
$s\geq 4$,
$\mathrel{\mathop{=\hskip-5pt=}\limits^{\hbox{or}}}(5,s)$,
$s\geq 5$,
$\mathrel{\mathop{=\hskip-5pt=}\limits^{\hbox{or}}}(6,s)$,
$s\geq 6$,
$\mathrel{\mathop{=\hskip-5pt=}\limits^{\hbox{or}}}(7,7)$.

\itemitem{(4)} $(a,b,c,d)=(3,2,s,r)$,
where $(s,r)=(u_{11},3)$, $u_{11}\in\{4,5,\ldots,11\}$,
$\mathrel{\mathop{=\hskip-5pt=}\limits^{\hbox{or}}}(u_{12},4)$,
$u_{12}\in\{5,6,7,8\}$,
$\mathrel{\mathop{=\hskip-5pt=}\limits^{\hbox{or}}}(u_{13},5)$,
$u_{13}\in\{6,7\}$,
$\mathrel{\mathop{=\hskip-5pt=}\limits^{\hbox{or}}}(7,6)$.

\itemitem{(5)} $(a,b,c,d)=(2,r,3,s)$,
where $(r,s)=(4,s)$, $s\geq 4$,
$\mathrel{\mathop{=\hskip-5pt=}\limits^{\hbox{or}}}(5,s)$,
$s\geq 5$,
$\mathrel{\mathop{=\hskip-5pt=}\limits^{\hbox{or}}}(6,s)$,
$s\geq 6$,
$\mathrel{\mathop{=\hskip-5pt=}\limits^{\hbox{or}}}(7,u_{14})$,
$u_{14}\in\{7,8,\ldots,18\}$,
$\mathrel{\mathop{=\hskip-5pt=}\limits^{\hbox{or}}}(8,u_{15})$,
$u_{15}\in\{8,9,10,11\}$.

\itemitem{(6)} $(a,b,c,d)=(2,r,s,3)$,
where $(r,s)=(4,u_{16})$, $u_{16}\in\{4,5,\ldots,15\}$,
$\mathrel{\mathop{=\hskip-5pt=}\limits^{\hbox{or}}}(5,u_{17})$,
$u_{17}\in\{5,6,\ldots,9\}$,
$\mathrel{\mathop{=\hskip-5pt=}\limits^{\hbox{or}}}(6,u_{18})$,
$u_{18}\in\{6,7\}$.

\itemitem{(7)} $(a,b,c,d)=(r,2,3,s)$,
where $(r,s)=(4,s)$, $s\geq 4$,
$\mathrel{\mathop{=\hskip-5pt=}\limits^{\hbox{or}}}(5,s)$,
$s\geq 5$,
$\mathrel{\mathop{=\hskip-5pt=}\limits^{\hbox{or}}}(6,s)$,
$s\geq 6$.

\itemitem{(8)} $(a,b,c,d)=(r,2,s,3)$,
where $(r,s)=(4,u_{19})$, $u_{19}\in\{4,5\}$.

\itemitem{(9)} $(a,b,c,d)=(2,s,3,r)$,
where $(s,r)=(s,3)$, $s\geq 4$,
$\mathrel{\mathop{=\hskip-5pt=}\limits^{\hbox{or}}}(s,4)$,
$s\geq 5$,
$\mathrel{\mathop{=\hskip-5pt=}\limits^{\hbox{or}}}(u_{20},5)$,
$u_{20}\in\{6,7,\ldots,20\}$,
$\mathrel{\mathop{=\hskip-5pt=}\limits^{\hbox{or}}}(u_{21},6)$,
$u_{21}\in\{7,8,\ldots,13\}$,
$\mathrel{\mathop{=\hskip-5pt=}\limits^{\hbox{or}}}(u_{22},7)$,
$u_{22}\in\{8,9,10\}$,
$\mathrel{\mathop{=\hskip-5pt=}\limits^{\hbox{or}}}(9,8)$.

\itemitem{(10)} $(a,b,c,d)=(2,s,r,3)$,
where $(s,r)=(s,4)$, $s\geq 5$,
$\mathrel{\mathop{=\hskip-5pt=}\limits^{\hbox{or}}}(u_{23},5)$,
$u_{23}\in\{6,7,\ldots,14\}$,
$\mathrel{\mathop{=\hskip-5pt=}\limits^{\hbox{or}}}(u_{24},6)$,
$u_{24}\in\{7,8\}$.

\itemitem{(11)} $(a,b,c,d)=(s,2,3,r)$,
where $(s,r)=(u_{25},3)$, $u_{25}\in\{4,5,6,7\}$,
$\mathrel{\mathop{=\hskip-5pt=}\limits^{\hbox{or}}}(u_{26},4)$,
$u_{26}\in\{5,6\}$,
$\mathrel{\mathop{=\hskip-5pt=}\limits^{\hbox{or}}}(6,5)$.

\itemitem{(12)} $(a,b,c,d)=(3,r,2,s)$,
where $(r,s)=(3,s)$, $s\geq 3$,
$\mathrel{\mathop{=\hskip-5pt=}\limits^{\hbox{or}}}(4,s)$,
$s\geq 4$,
$\mathrel{\mathop{=\hskip-5pt=}\limits^{\hbox{or}}}(5,s)$,
$s\geq 5$,
$\mathrel{\mathop{=\hskip-5pt=}\limits^{\hbox{or}}}(6,s)$,
$s\geq 6$,
$\mathrel{\mathop{=\hskip-5pt=}\limits^{\hbox{or}}}(7,u_{27})$,
$u_{27}\in\{7,8,\ldots,12\}$.

\itemitem{(13)} $(a,b,c,d)=(3,r,s,2)$,
where $(r,s)=(3,u_{28})$, $u_{28}\in\{3,4,\ldots,8\}$,
$\mathrel{\mathop{=\hskip-5pt=}\limits^{\hbox{or}}}(4,u_{29})$,
$u_{29}\in\{4,5\}$.

\itemitem{(14)} $(a,b,c,d)=(r,3,2,s)$,
where $(r,s)=(4,s)$, $s\geq 4$,
$\mathrel{\mathop{=\hskip-5pt=}\limits^{\hbox{or}}}(5,s)$,
$s\geq 5$,
$\mathrel{\mathop{=\hskip-5pt=}\limits^{\hbox{or}}}(6,s)$,
$s\geq 6$.

\itemitem{(15)} $(a,b,c,d)=(3,s,2,r)$,
where $(s,r)=(s,3)$, $s\geq 4$,
$\mathrel{\mathop{=\hskip-5pt=}\limits^{\hbox{or}}}(u_{30},4)$,
$u_{30}\in\{5,6,\ldots,15\}$,
$\mathrel{\mathop{=\hskip-5pt=}\limits^{\hbox{or}}}(u_{31},5)$,
$u_{31}\in\{6,7,\ldots,10\}$,
$\mathrel{\mathop{=\hskip-5pt=}\limits^{\hbox{or}}}(u_{32},6)$,
$u_{32}\in\{7,8,9\}$,
$\mathrel{\mathop{=\hskip-5pt=}\limits^{\hbox{or}}}(8,7)$.

\itemitem{(16)} $(a,b,c,d)=(3,s,r,2)$,
where $(s,r)=(s,3)$, $s\geq 4$,
$\mathrel{\mathop{=\hskip-5pt=}\limits^{\hbox{or}}}(u_{33},4)$,
$u_{33}\in\{5,6,7\}$.

\itemitem{(17)} $(a,b,c,d)=(s,3,2,r)$,
where $(s,r)=(u_{34},3)$, $u_{34}\in\{4,5,6,7\}$,
$\mathrel{\mathop{=\hskip-5pt=}\limits^{\hbox{or}}}(u_{35},4)$,
$u_{35}\in\{5,6,7\}$,
$\mathrel{\mathop{=\hskip-5pt=}\limits^{\hbox{or}}}(6,5)$.

\itemitem{(18)} $(a,b,c,d)=(4,3,3,2)$.

\itemitem{(19)} $(a,b,c,d)=(r,s,2,3)$,
where $(r,s)=(4,u_{36})$, $u_{36}\in\{4,5,6\}$.

\itemitem{(20)} $(a,b,c,d)=(4,4,3,2)$.

\itemitem{(21)} $(a,b,c,d)=(5,4,2,3)$.

\smallskip

\item{4.}
\Itemitem{(1)} $(a,b,c,d)=(2,4,r,s)$,
where $(r,s)=(4,s)$, $s\geq 4$,
$\mathrel{\mathop{=\hskip-5pt=}\limits^{\hbox{or}}}(5,u_{37})$,
$u_{37}\in\{5,6,7,8\}$.

\itemitem{(2)} $(a,b,c,d)=(2,4,s,r)$,
where $(s,r)=(u_{38},4)$, $u_{38}\in\{5,6,7\}$,
$\mathrel{\mathop{=\hskip-5pt=}\limits^{\hbox{or}}}(6,5)$.

\itemitem{(3)} $(a,b,c,d)=(4,2,r,s)$,
where $(r,s)=(4,s)$, $s\geq 4$.

\itemitem{(4)} $(a,b,c,d)=(2,r,4,s)$,
where $(r,s)=(5,u_{39})$, $u_{39}\in\{5,6,7,8\}$.

\itemitem{(5)} $(a,b,c,d)=(2,5,5,4)$.

\itemitem{(6)} $(a,b,c,d)=(2,s,4,r)$,
where $(s,r)=(u_{40},4)$, $u_{40}\in\{5,6,\ldots,9\}$,
$\mathrel{\mathop{=\hskip-5pt=}\limits^{\hbox{or}}}(6,5)$.

\itemitem{(7)} $(a,b,c,d)=(2,6,5,4)$.

\itemitem{(8)} $(a,b,c,d)=(4,r,2,s)$,
where $(r,s)=(4,s)$, $s\geq 4$.

\itemitem{(9)} $(a,b,c,d)=(4,5,2,4)$.

\smallskip

\item{5.}
\Itemitem{(1)} $(a,b,c,d)=(3,3,r,s)$,
where $(r,s)=(3,s)$, $s\geq 3$.

\itemitem{(2)} $(a,b,c,d)=(3,3,4,3)$.

\itemitem{(3)} $(a,b,c,d)=(3,4,3,3)$.

\medskip

\noindent (X) $f(x,y,z,w)=x^a+y^bz+z^cw+yw^d$.

\smallskip
\item{1.}
\Itemitem{(1)}
$(a,b,c,d)=(r,1,s,t)$, $r\geq 2$, $s\geq 2$, $t\geq 2$.

\itemitem{(2)} $(a,b,c,d)=(r,1,1,s)$, $r\geq 2$, $s\geq 2$.

\itemitem{(3)} $(a,b,c,d)=(r,1,1,1)$, $r\geq 2$.

\smallskip

\item{2.}
\Itemitem{(1)} $a=2$, $(b,c,d)=(2,r,s)$, $r\geq 2$, $s\geq 2$.
\itemitem{(2)} $a=2$, $(b,c,d)=(2,s,r)$, $s\geq r+1$, $r\geq 3$.
\itemitem{(3)} $a=r$, $(b,c,d)=(2,2,s)$,
where $(r,s)=(3,s)$, $s\geq 3$,
$\mathrel{\mathop{=\hskip-5pt=}\limits^{\hbox{or}}}(4,s)$,
$s\geq 4$,
$\mathrel{\mathop{=\hskip-5pt=}\limits^{\hbox{or}}}(5,u_1)$,
$u_1\in\{5,6,\ldots,10\}$,
$\mathrel{\mathop{=\hskip-5pt=}\limits^{\hbox{or}}}(6,6)$.

\itemitem{(4)} $a=s$, $(b,c,d)=(2,2,r)$,
where $(s,r)=(s,2)$, $s\geq 3$,
$\mathrel{\mathop{=\hskip-5pt=}\limits^{\hbox{or}}}(u_2,3)$,
$u_2\in\{4,5,\ldots,12\}$,
$\mathrel{\mathop{=\hskip-5pt=}\limits^{\hbox{or}}}(u_3,4)$,
$u_3\in\{5,6,7,8\}$,
$\mathrel{\mathop{=\hskip-5pt=}\limits^{\hbox{or}}}(6,5)$.

\smallskip

\item{3.}
\Itemitem{(1)} $a=2$, $(b,c,d)=(3,r,s)$,
where $(r,s)=(3,s)$, $s\geq 3$,
$\mathrel{\mathop{=\hskip-5pt=}\limits^{\hbox{or}}}(4,s)$,
$s\geq 4$,
$\mathrel{\mathop{=\hskip-5pt=}\limits^{\hbox{or}}}(5,u_7)$,
$u_7\in\{5,6,\ldots,18\}$,
$\mathrel{\mathop{=\hskip-5pt=}\limits^{\hbox{or}}}(6,u_8)$,
$u_8\in\{6,7,\ldots,11\}$,
$\mathrel{\mathop{=\hskip-5pt=}\limits^{\hbox{or}}}(7,u_9)$,
$u_9\in\{7,8\}$.

\itemitem{(2)} $a=2$, $(b,c,d)=(3,s,r)$,
where $(s,r)=(s,4)$, $s\geq 5$,
$\mathrel{\mathop{=\hskip-5pt=}\limits^{\hbox{or}}}(u_{10},5)$,
$u_{10}\in\{6,7,\ldots,18\}$,
$\mathrel{\mathop{=\hskip-5pt=}\limits^{\hbox{or}}}(u_{11},6)$,
$u_{11}\in\{7,8,\ldots,11\}$,
$\mathrel{\mathop{=\hskip-5pt=}\limits^{\hbox{or}}}(8,7)$.

\itemitem{(3)} $a=3$, $(b,c,d)=(2,r,s)$,
where $(r,s)=(3,s)$, $s\geq 3$,
$\mathrel{\mathop{=\hskip-5pt=}\limits^{\hbox{or}}}(4,u_{12})$,
$u_{12}\in\{4,5,\ldots,12\}$,
$\mathrel{\mathop{=\hskip-5pt=}\limits^{\hbox{or}}}(5,u_{13})$,
$u_{13}\in\{5,6,7\}$,
$\mathrel{\mathop{=\hskip-5pt=}\limits^{\hbox{or}}}(6,6)$.

\itemitem{(4)} $a=3$, $(b,c,d)=(2,s,r)$,
where $(s,r)=(s,3)$, $s\geq 4$,
$\mathrel{\mathop{=\hskip-5pt=}\limits^{\hbox{or}}}(u_{14},4)$,
$u_{14}\in\{5,6,\ldots,12\}$,
$\mathrel{\mathop{=\hskip-5pt=}\limits^{\hbox{or}}}(u_{15},5)$,
$u_{15}\in\{6,7\}$.

\itemitem{(5)} $a=r$, $(b,c,d)=(2,3,s)$,
where $(r,s)=(4,u_{16})$, $u_{16}\in\{4,5,6\}$.

\itemitem{(6)} $a=r$, $(b,c,d)=(2,s,3)$,
where $(r,s)=(4,u_{17})$, $u_{17}\in\{4,5,6\}$.

\itemitem{(7)} $a=s$, $(b,c,d)=(2,3,r)$,
where $(s,r)=(u_{18},3)$, $u_{18}\in\{4,5,6\}$.

\smallskip
\item{4.}
\Itemitem{(1)} $a=2$, $(b,c,d)=(4,r,s)$,
where $(r,s)=(4,u_{18})$, $u_{18}\in\{4,5,\ldots,10\}$,
$\mathrel{\mathop{=\hskip-5pt=}\limits^{\hbox{or}}}(5,u_{19})$,
$u_{19}\in\{5,6\}$.

\itemitem{(2)} $a=2$, $(b,c,d)=(4,6,5)$.

\itemitem{(3)} $a=4$, $(b,c,d)=(2,4,4)$.

\smallskip

\item{5.}
$a=3$, $(b,c,d)=(3,r,s)$,
where $(r,s)=(3,u_{20})$, $u_{20}\in\{3,4,5\}$.

\medskip

\noindent (XI) $f(x,y,z,w)=x^a+xy^b+yz^c+zw^d$.

\smallskip
\item{1.}
\Itemitem{(1)}
$(a,b,c,d)=(r,1,s,t)$, $r\geq 2$, $s\geq 2$, $t\geq 2$.
\itemitem{(2)} $(a,b,c,d)=(r,s,1,t)$, $r\geq 2$, $s\geq 2$, $t\geq 2$.
\itemitem{(3)} $(a,b,c,d)=(r,s,t,1)$, $r\geq 2$, $s\geq 2$, $t\geq 2$.
\itemitem{(4)} $(a,b,c,d)=(r,1,1,t)$, $r\geq 2$, $s\geq 2$.
\itemitem{(5)} $(a,b,c,d)=(r,1,s,1)$, $r\geq 2$, $s\geq 2$.
\itemitem{(6)} $(a,b,c,d)=(r,s,1,1)$, $r\geq 2$, $s\geq 2$.
\itemitem{(7)} $(a,b,c,d)=(r,1,1,1)$, $r\geq 2$.

\smallskip
\item{2.}
\Itemitem{(1)} $(a,b,c,d)=(2,2,r,s)$,
where $(r,s)=(2,s)$, $s\geq 2$,
$\mathrel{\mathop{=\hskip-5pt=}\limits^{\hbox{or}}}(3,s)$,
$s\geq 3$,
$\mathrel{\mathop{=\hskip-5pt=}\limits^{\hbox{or}}}(4,u_1)$,
$u_1\in\{4,5,\ldots,12\}$,
$\mathrel{\mathop{=\hskip-5pt=}\limits^{\hbox{or}}}(5,u_2)$,
$u_2\in\{5,6,7,8\}$,
$\mathrel{\mathop{=\hskip-5pt=}\limits^{\hbox{or}}}(6,6)$.

\itemitem{(2)} $(a,b,c,d)=(2,2,s,r)$,
where $(s,r)=(s,2)$, $s\geq 3$,
$\mathrel{\mathop{=\hskip-5pt=}\limits^{\hbox{or}}}(s,3)$,
$s\geq 4$,
$\mathrel{\mathop{=\hskip-5pt=}\limits^{\hbox{or}}}(s,4)$,
$s\geq 5$,
$\mathrel{\mathop{=\hskip-5pt=}\limits^{\hbox{or}}}(u_3,5)$,
$u_3\in\{6,7,\ldots,11\}$,
$\mathrel{\mathop{=\hskip-5pt=}\limits^{\hbox{or}}}(7,6)$.

\itemitem{(3)} $(a,b)=(2,r)$, $\{c,d\}=\{2,s\}$, $r\geq 3$, $s\geq r$.

\itemitem{(4)} $(a,b,c,d)=(r,2,2,s)$,
where $(r,s)=(3,s)$, $s\geq 3$,
$\mathrel{\mathop{=\hskip-5pt=}\limits^{\hbox{or}}}(4,u_4)$,
$u_4\in\{4,5,\ldots,10\}$,
$\mathrel{\mathop{=\hskip-5pt=}\limits^{\hbox{or}}}(5,u_5)$,
$u_5\in\{5,6\}$.

\itemitem{(5)} $(a,b,c,d)=(r,2,s,2)$, $r\geq 3$, $s\geq r$.

\itemitem{(6)} $(a,b)=(2,s)$, $\{c,d\}=\{2,r\}$, $r\geq 2$, $s\geq r+1$.

\itemitem{(7)} $(a,b,c,d)=(s,2,2,r)$,
where $(s,r)=(s,2)$, $s\geq 3$,
$\mathrel{\mathop{=\hskip-5pt=}\limits^{\hbox{or}}}(s,3)$,
$s\geq 4$,
$\mathrel{\mathop{=\hskip-5pt=}\limits^{\hbox{or}}}(u_6,4)$,
$u_6\in\{5,6,\ldots,10\}$,
$\mathrel{\mathop{=\hskip-5pt=}\limits^{\hbox{or}}}(6,5)$.

\itemitem{(8)} $(a,b,c,d)=(s,2,r,2)$, $r\geq 3$, $s\geq r+1$.

\itemitem{(9)} $(a,b,c,d)=(r,s,2,2)$,
where $(r,s)=(3,s)$, $s\geq 3$,
$\mathrel{\mathop{=\hskip-5pt=}\limits^{\hbox{or}}}(4,s)$,
$s\geq 4$,
$\mathrel{\mathop{=\hskip-5pt=}\limits^{\hbox{or}}}(5,u_7)$,
$u_7\in\{5,6,\ldots,11\}$,
$\mathrel{\mathop{=\hskip-5pt=}\limits^{\hbox{or}}}(6,u_8)$,
$u_8\in\{6,7\}$.

\itemitem{(10)} $(a,b,c,d)=(s,r,2,2)$,
where $(s,r)=(s,3)$, $s\geq 4$,
$\mathrel{\mathop{=\hskip-5pt=}\limits^{\hbox{or}}}(u_9,4)$,
$u_9\in\{5,6,\ldots,12\}$,
$\mathrel{\mathop{=\hskip-5pt=}\limits^{\hbox{or}}}(u_{10},5)$,
$u_{10}\in\{6,7,8\}$.

\smallskip
\item{3.}
\Itemitem{(1)} $(a,b,c,d)=(2,3,r,s)$,
where $(r,s)=(3,u_{11})$, $u_{11}\in\{3,4,\ldots,12\}$,
$\mathrel{\mathop{=\hskip-5pt=}\limits^{\hbox{or}}}(4,u_{12})$,
$u_{12}\in\{4,5,6\}$.

\itemitem{(2)} $(a,b,c,d)=(2,3,s,r)$,
where $(s,r)=(s,3)$, $s\geq 4$,
$\mathrel{\mathop{=\hskip-5pt=}\limits^{\hbox{or}}}(u_{13},4)$,
$u_{13}\in\{5,6,7\}$.

\itemitem{(3)} $(a,b,c,d)=(3,2,r,s)$,
where $(r,s)=(3,u_{14})$, $u_{14}\in\{3,4,5,6\}$,
$\mathrel{\mathop{=\hskip-5pt=}\limits^{\hbox{or}}}(4,4)$.

\itemitem{(4)} $(a,b,c,d)=(3,2,s,r)$,
where $(s,r)=(s,3)$, $s\geq 4$,
$\mathrel{\mathop{=\hskip-5pt=}\limits^{\hbox{or}}}(5,4)$.

\itemitem{(5)} $(a,b,c,d)=(2,r,3,s)$,
where $(r,s)=(4,u_{15})$, $u_{15}\in\{4,5,\ldots,8\}$,
$\mathrel{\mathop{=\hskip-5pt=}\limits^{\hbox{or}}}(5,u_{16})$,
$u_{16}\in\{5,6\}$,
$\mathrel{\mathop{=\hskip-5pt=}\limits^{\hbox{or}}}(6,6)$.

\itemitem{(6)} $(a,b,c,d)=(2,r,s,3)$,
where $(r,s)=(4,u_{17})$, $u_{17}\in\{4,5,\ldots,12\}$,
$\mathrel{\mathop{=\hskip-5pt=}\limits^{\hbox{or}}}(5,u_{18})$,
$u_{18}\in\{5,6,7,8\}$,
$\mathrel{\mathop{=\hskip-5pt=}\limits^{\hbox{or}}}(6,u_{19})$,
$u_{19}\in\{6,7\}$.

\itemitem{(7)} $(a,b,c,d)=(r,2,3,s)$,
where $(r,s)=(4,4)$.

\itemitem{(8)} $(a,b,c,d)=(r,2,s,3)$,
where $(r,s)=(4,u_{20})$, $u_{20}\in\{4,5,\ldots,9\}$,
$\mathrel{\mathop{=\hskip-5pt=}\limits^{\hbox{or}}}(5,5)$.

\itemitem{(9)} $(a,b,c,d)=(2,s,3,r)$,
where $(s,r)=(s,3)$, $s\geq 4$,
$\mathrel{\mathop{=\hskip-5pt=}\limits^{\hbox{or}}}(s,4)$,
$s\geq 5$,
$\mathrel{\mathop{=\hskip-5pt=}\limits^{\hbox{or}}}(u_{21},5)$,
$u_{21}\in\{6,7,\ldots,10\}$.

\itemitem{(10)} $(a,b,c,d)=(2,s,r,3)$,
where $(s,r)=(s,4)$, $s\geq 5$,
$\mathrel{\mathop{=\hskip-5pt=}\limits^{\hbox{or}}}(u_{22},5)$,
$u_{22}\in\{6,7,\ldots,12\}$,
$\mathrel{\mathop{=\hskip-5pt=}\limits^{\hbox{or}}}(7,6)$.

\itemitem{(11)} $(a,b,c,d)=(s,2,3,r)$,
where $(s,r)=(u_{23},3)$, $u_{23}\in\{4,5,\ldots,10\}$.

\itemitem{(12)} $(a,b,c,d)=(s,2,r,3)$,
where $(s,r)=(u_{24},4)$, $u_{24}\in\{5,6\}$.

\itemitem{(13)} $(a,b,c,d)=(3,r,2,s)$,
where $(r,s)=(3,u_{25})$, $u_{25}\in\{3,4,\ldots,10\}$,
$\mathrel{\mathop{=\hskip-5pt=}\limits^{\hbox{or}}}(4,u_{26})$,
$u_{26}\in\{4,5,6\}$,
$\mathrel{\mathop{=\hskip-5pt=}\limits^{\hbox{or}}}(5,5)$.

\itemitem{(14)} $(a,b,c,d)=(3,r,s,2)$,
where $(r,s)=(3,s)$, $s\geq 3$,
$\mathrel{\mathop{=\hskip-5pt=}\limits^{\hbox{or}}}(4,s)$,
$s\geq 4$,
$\mathrel{\mathop{=\hskip-5pt=}\limits^{\hbox{or}}}(5,u_{27})$,
$u_{27}\in\{5,6,\ldots,12\}$,
$\mathrel{\mathop{=\hskip-5pt=}\limits^{\hbox{or}}}(6,u_{28})$,
$u_{28}\in\{6,7\}$.

\itemitem{(15)} $(a,b,c,d)=(4,3,2,4)$.

\itemitem{(16)} $(a,b,c,d)=(r,3,s,2)$,
where $(r,s)=(4,s)$, $s\geq 4$,
$\mathrel{\mathop{=\hskip-5pt=}\limits^{\hbox{or}}}(5,u_{29})$,
$u_{29}\in\{5,6,\ldots,10\}$,
$\mathrel{\mathop{=\hskip-5pt=}\limits^{\hbox{or}}}(6,6)$.

\itemitem{(17)} $(a,b,c,d)=(3,s,2,r)$,
where $(s,r)=(s,3)$, $s\geq 4$,
$\mathrel{\mathop{=\hskip-5pt=}\limits^{\hbox{or}}}(u_{30},4)$,
$u_{30}\in\{5,6,\ldots,9\}$.

\itemitem{(18)} $(a,b,c,d)=(3,s,r,2)$,
where $(s,r)=(s,3)$, $s\geq 4$,
$\mathrel{\mathop{=\hskip-5pt=}\limits^{\hbox{or}}}(u_{31},4)$,
$u_{31}\in\{5,6,\ldots,13\}$,
$\mathrel{\mathop{=\hskip-5pt=}\limits^{\hbox{or}}}(u_{32},5)$,
$u_{32}\in\{6,7,8\}$,
$\mathrel{\mathop{=\hskip-5pt=}\limits^{\hbox{or}}}(7,6)$.

\itemitem{(19)} $(a,b,c,d)=(s,3,2,r)$,
where $(s,r)=(u_{33},3)$, $u_{33}\in\{4,5,6\}$.

\itemitem{(20)} $(a,b,c,d)=(s,3,r,2)$,
where $(s,r)=(u_{34},3)$, $u_{34}\in\{4,5,\ldots,12\}$,
$\mathrel{\mathop{=\hskip-5pt=}\limits^{\hbox{or}}}(u_{35},4)$,
$u_{35}\in\{5,6,7,8\}$,
$\mathrel{\mathop{=\hskip-5pt=}\limits^{\hbox{or}}}(6,5)$.

\itemitem{(21)} $(a,b,c,d)=(r,s,2,3)$,
where $(r,s)=(4,u_{36})$, $u_{36}\in\{4,5\}$.

\itemitem{(22)} $(a,b,c,d)=(r,s,3,2)$,
where $(r,s)=(4,u_{37})$, $u_{37}\in\{4,5,6,7\}$.

\itemitem{(23)} $(a,b,c,d)=(s,r,3,2)$,
where $(s,r)=(u_{38},4)$, $u_{38}\in\{5,6\}$.

\smallskip
\item{4.}
\Itemitem{(1)} $(a,b,c,d)=(2,4,4,4)$.

\itemitem{(2)} $(a,b,c,d)=(2,4,5,4)$.

\itemitem{(3)} $(a,b,c,d)=(2,u_{39},4,4)$, $u_{39}\in\{5,6\}$.

\itemitem{(4)} $(a,b,c,d)=(4,4,u_{40},2)$, $u_{40}\in\{4,5,6\}$.

\itemitem{(5)} $(a,b,c,d)=(4,5,4,2)$.

\smallskip
\item{5.}
\Itemitem{(1)} $(a,b,c,d)=(3,3,3,3)$.

\itemitem{(2)} $(a,b,c,d)=(3,3,4,3)$.

\itemitem{(3)} $(a,b,c,d)=(3,4,3,3)$.

\medskip

\noindent (XII) $f(x,y,z,w)=x^a+xy^b+xz^c+yw^d+y^pz^q$,
$\frac{p(a-1)}{ab}+\frac{q(a-1)}{ac}=1$.

\smallskip
\item{1.}
\Itemitem{(1)}
$(a,b,c,d)=(r,1,s,t)$, $r\geq 2$, $s\geq 2$, $t\geq 2$.
\itemitem{(2)} $(a,b,c,d)=(r,s,1,t)$, $r\geq 2$, $s\geq 2$, $t\geq 2$.
\itemitem{(3)} $(a,b,c,d)=(r,s,t,1)$, $r\geq 2$, $s\geq 2$, $t\geq 2$.
\itemitem{(4)} $(a,b,c,d)=(r,1,1,s)$, $r\geq 2$, $s\geq 2$.
\itemitem{(5)} $(a,b,c,d)=(r,1,s,1)$, $r\geq 2$, $s\geq 2$.
\itemitem{(6)} $(a,b,c,d)=(r,s,1,1)$, $r\geq 2$, $s\geq 2$.
\itemitem{(7)} $(a,b,c,d)=(r,1,1,1)$, $r\geq 2$.

\smallskip
\item{2.}
\Itemitem{(1)} $(a,b,c,d)=(2,2,r,s)$,
where $(r,s)=(2,s)$, $s\geq 2$,
$\mathrel{\mathop{=\hskip-5pt=}\limits^{\hbox{or}}}(3,u_1)$,
$u_1\in\{3,4,\ldots,8\}$,
$\mathrel{\mathop{=\hskip-5pt=}\limits^{\hbox{or}}}(4,u_2)$,
$u_2\in\{4,5\}$.

\itemitem{(2)} $(a,b,c,d)=(2,2,s,r)$,
where $(s,r)=(s,2)$, $s\geq 3$,
$\mathrel{\mathop{=\hskip-5pt=}\limits^{\hbox{or}}}(s,3)$,
$s\geq 4$,
$\mathrel{\mathop{=\hskip-5pt=}\limits^{\hbox{or}}}(u_3,4)$,
$u_3\in\{5,6,7\}$.

\itemitem{(3)} $(a,b,c,d)=(2,r,2,s)$,
where $(r,s)=(3,u_4)$, $u_4\in\{3,4,\ldots,9\}$,
$\mathrel{\mathop{=\hskip-5pt=}\limits^{\hbox{or}}}(4,u_5)$,
$u_5\in\{4,5,6\}$,
$\mathrel{\mathop{=\hskip-5pt=}\limits^{\hbox{or}}}(5,5)$.

\itemitem{(4)} $(a,b,c,d)=(2,r,s,2)$, $r\geq 3$, $s\geq r$.

\itemitem{(5)} $(a,b,c,d)=(r,2,2,s)$, $r\geq 3$, $s\geq r$.

\itemitem{(6)} $(a,b,c,d)=(r,2,s,2)$,
where $(r,s)=(3,s)$, $s\geq 3$,
$\mathrel{\mathop{=\hskip-5pt=}\limits^{\hbox{or}}}(4,u_6)$,
$u_6\in\{4,5,\ldots,11\}$,
$\mathrel{\mathop{=\hskip-5pt=}\limits^{\hbox{or}}}(5,u_7)$,
$u_7\in\{5,6,7\}$,
$\mathrel{\mathop{=\hskip-5pt=}\limits^{\hbox{or}}}(6,6)$.

\itemitem{(7)} $(a,b,c,d)=(2,s,2,r)$,
where $(s,r)=(s,2)$, $s\geq 3$,
$\mathrel{\mathop{=\hskip-5pt=}\limits^{\hbox{or}}}(s,3)$,
$s\geq 4$,
$\mathrel{\mathop{=\hskip-5pt=}\limits^{\hbox{or}}}(s,4)$,
$s\geq 5$,
$\mathrel{\mathop{=\hskip-5pt=}\limits^{\hbox{or}}}(u_8,5)$,
$u_8\in\{6,7\}$.

\itemitem{(8)} $(a,b,c,d)=(2,s,r,2)$, $r\geq 3$, $s\geq r+1$.

\itemitem{(9)} $(a,b,c,d)=(s,2,2,r)$, $r\geq 2$, $s\geq r+1$.

\itemitem{(10)} $(a,b,c,d)=(s,2,r,2)$,
where $(s,r)=(s,3)$, $s\geq 4$,
$\mathrel{\mathop{=\hskip-5pt=}\limits^{\hbox{or}}}(s,4)$,
$s\geq 5$,
$\mathrel{\mathop{=\hskip-5pt=}\limits^{\hbox{or}}}(u_{9},5)$,
$u_{9}\in\{6,7,\ldots,10\}$.

\itemitem{(11)} $\{a,b\}=\{r,s\}$, $\{c,d\}=\{2,2\}$, $r\geq 3$, $s\geq r$.

\smallskip
\item{3.}
\Itemitem{(1)} $(a,b,c,d)=(2,3,r,s)$,
where $(r,s)=(3,u_{10})$, $u_{10}\in\{3,4\}$.

\itemitem{(2)} $(a,b,c,d)=(2,3,s,r)$,
where $(s,r)=(u_{11},3)$, $u_{11}\in\{4,5,\ldots,8\}$.

\itemitem{(3)} $(a,b,c,d)=(3,2,r,s)$,
where $(r,s)=(3,u_{12})$, $u_{12}\in\{3,4,5\}$.

\itemitem{(4)} $(a,b,c,d)=(3,2,s,r)$,
where $(s,r)=(u_{13},3)$, $u_{13}\in\{4,5\}$.

\itemitem{(5)} $(a,b,c,d)=(2,4,3,4)$.

\itemitem{(6)} $(a,b,c,d)=(2,r,s,3)$,
where $(r,s)=(4,u_{14})$, $u_{14}\in\{4,5\}$.

\itemitem{(7)} $(a,b,c,d)=(4,2,3,4)$.

\itemitem{(8)} $(a,b,c,d)=(4,2,4,3)$.

\itemitem{(9)} $(a,b,c,d)=(2,s,3,r)$,
where $(s,r)=(s,3)$, $s\geq 4$.

\itemitem{(10)} $(a,b,c,d)=(2,s,r,3)$,
where $(s,r)=(u_{15},3)$, $u_{15}\in\{5,6,7\}$.

\itemitem{(11)} $(a,b,c,d)=(s,2,3,r)$,
where $(s,r)=(s,3)$, $s\geq 4$,
$\mathrel{\mathop{=\hskip-5pt=}\limits^{\hbox{or}}}(u_{16},4)$,
$u_{16}\in\{5,6\}$.

\itemitem{(12)} $(a,b,c,d)=(3,r,2,s)$,
where $(r,s)=(3,u_{17})$, $u_{17}\in\{3,4,5,6\}$,
$\mathrel{\mathop{=\hskip-5pt=}\limits^{\hbox{or}}}(4,4)$.

\itemitem{(13)} $(a,b,c,d)=(3,r,s,2)$,
where $(r,s)=(3,u_{18})$, $u_{18}\in\{3,4,\ldots,11\}$,
$\mathrel{\mathop{=\hskip-5pt=}\limits^{\hbox{or}}}(4,u_{19})$,
$u_{19}\in\{4,5,6,7\}$,
$\mathrel{\mathop{=\hskip-5pt=}\limits^{\hbox{or}}}(5,u_{20})$,
$u_{20}\in\{5,6\}$.

\itemitem{(14)} $(a,b,c,d)=(r,3,2,s)$,
where $(r,s)=(4,u_{21})$, $u_{21}\in\{4,5\}$,
$\mathrel{\mathop{=\hskip-5pt=}\limits^{\hbox{or}}}(5,5)$.

\itemitem{(15)} $(a,b,c,d)=(r,3,s,2)$,
where $(r,s)=(4,u_{22})$, $u_{22}\in\{4,5\}$.

\itemitem{(16)} $(a,b,c,d)=(3,s,2,r)$,
where $(s,r)=(s,3)$, $s\geq 4$,
$\mathrel{\mathop{=\hskip-5pt=}\limits^{\hbox{or}}}(4,5)$.

\itemitem{(17)} $(a,b,c,d)=(3,s,r,2)$,
where $(s,r)=(s,3)$, $s\geq 4$,
$\mathrel{\mathop{=\hskip-5pt=}\limits^{\hbox{or}}}(s,4)$,
$s\geq 5$,
$\mathrel{\mathop{=\hskip-5pt=}\limits^{\hbox{or}}}(u_{23},5)$,
$u_{23}\in\{6,7,8,9\}$.

\itemitem{(18)} $(a,b,c,d)=(s,3,2,r)$,
where $(s,r)=(s,3)$, $s\geq 4$,
$\mathrel{\mathop{=\hskip-5pt=}\limits^{\hbox{or}}}(s,4)$,
$s\geq 5$,
$\mathrel{\mathop{=\hskip-5pt=}\limits^{\hbox{or}}}(6,5)$.

\itemitem{(19)} $(a,b,c,d)=(s,3,r,2)$,
where $(s,r)=(s,3)$, $s\geq 4$,
$\mathrel{\mathop{=\hskip-5pt=}\limits^{\hbox{or}}}(u_{24},4)$,
$u_{24}\in\{5,6\}$.

\itemitem{(20)} $(a,b,c,d)=(r,s,2,3)$,
where $(r,s)=(4,u_{25})$, $u_{25}\in\{4,5,\ldots,11\}$,
$\mathrel{\mathop{=\hskip-5pt=}\limits^{\hbox{or}}}(5,u_{26})$,
$u_{26}\in\{5,6,7\}$,
$\mathrel{\mathop{=\hskip-5pt=}\limits^{\hbox{or}}}(6,6)$.

\itemitem{(21)} $(a,b,c,d)=(r,s,3,2)$,
where $(r,s)=(4,s)$, $s\geq 4$,
$\mathrel{\mathop{=\hskip-5pt=}\limits^{\hbox{or}}}(5,u_{27})$,
$u_{27}\in\{5,6,\ldots,11\}$,
$\mathrel{\mathop{=\hskip-5pt=}\limits^{\hbox{or}}}(6,u_{28})$,
$u_{28}\in\{6,7\}$.

\itemitem{(22)} $(a,b,c,d)=(s,r,2,3)$,
where $(s,r)=(s,4)$, $s\geq 5$,
$\mathrel{\mathop{=\hskip-5pt=}\limits^{\hbox{or}}}(u_{30},5)$,
$u_{30}\in\{6,7,8\}$.

\smallskip
\item{4.}
\Itemitem{(1)} $(a,b,c,d)=(4,4,2,4)$.

\itemitem{(2)} $(a,b,c,d)=(4,4,4,2)$.

\itemitem{(3)} $(a,b,c,d)=(4,5,4,2)$.

\smallskip
\item{5.}
$(a,b,c,d)=(3,3,3,3)$.

\medskip

\noindent (XIII) $f(x,y,z,w)=x^a+xy^b+yz^c+yw^d+z^pw^q$,
$\frac{p(a(b-1)+1)}{abc}+\frac{q(a(b-1)+1}{abd}=1$.

\smallskip
\item{1.}
\Itemitem{(1)}
$(a,b,c,d)=(r,1,s,t)$, $r\geq 2$, $s\geq 2$, $t\geq 2$.
\itemitem{(2)} $(a,b,c,d)=(r,s,1,t)$, $r\geq 2$, $s\geq 2$, $t\geq 2$.
\itemitem{(3)} $(a,b,c,d)=(r,s,t,1)$, $r\geq 2$, $s\geq 2$, $t\geq 2$.
\itemitem{(4)} $(a,b,c,d)=(r,1,1,s)$, $r\geq 2$, $s\geq 2$.
\itemitem{(5)} $(a,b,c,d)=(r,1,s,1)$, $r\geq 2$, $s\geq 2$.
\itemitem{(6)} $(a,b,c,d)=(r,s,1,1)$, $r\geq 2$, $s\geq 2$.
\itemitem{(7)} $(a,b,c,d)=(r,1,1,1)$, $r\geq 2$.

\smallskip
\item{2.}
\Itemitem{(1)} $(a,b)=(2,2)$, $\{c,d\}=\{r,s\}$,
where $(r,s)=(2,s)$, $s\geq 2$,
$\mathrel{\mathop{=\hskip-5pt=}\limits^{\hbox{or}}}(3,s)$,
$s\geq 3$,
$\mathrel{\mathop{=\hskip-5pt=}\limits^{\hbox{or}}}(4,u_1)$,
$u_1\in\{4,5,\ldots,11\}$,
$\mathrel{\mathop{=\hskip-5pt=}\limits^{\hbox{or}}}(5,u_2)$,
$u_2\in\{5,6,7\}$.

\itemitem{(2)} $(a,b)=(2,r)$, $\{c,d\}=\{2,s\}$, $r\geq 3$, $s\geq r$.

\itemitem{(3)} $(a,b)=(r,2)$, $\{c,d\}=\{2,s\}$,
where $(r,s)=(3,s)$, $s\geq 3$,
$\mathrel{\mathop{=\hskip-5pt=}\limits^{\hbox{or}}}(4,u_3)$,
$u_3\in\{4,5,\ldots,9\}$,
$\mathrel{\mathop{=\hskip-5pt=}\limits^{\hbox{or}}}(5,5)$.

\itemitem{(4)} $(a,b)=(2,s)$, $\{c,d\}=\{2,r\}$, $r\geq 2$, $s\geq r+1$.

\itemitem{(5)} $(a,b)=(s,2)$, $\{c,d\}=\{2,r\}$,
where $(s,r)=\{s,2\}$, $s\geq 3$,
$\mathrel{\mathop{=\hskip-5pt=}\limits^{\hbox{or}}}(u_4,3)$,
$u_4\in\{4,5,\ldots,10\}$,
$\mathrel{\mathop{=\hskip-5pt=}\limits^{\hbox{or}}}(u_5,4)$,
$u_5\in\{5,6\}$.

\itemitem{(6)} $(a,b)=\{r,s\}$, $\{c,d\}=\{2,2\}$, $r\geq 3$, $s\geq r$.

\smallskip
\item{3.}
\Itemitem{(1)} $(a,b)=(2,3)$, $\{c,d\}=\{r,s\}$,
where $(r,s)=(3,u_6)$, $u_6\in\{3,4,\ldots,14\}$,
$\mathrel{\mathop{=\hskip-5pt=}\limits^{\hbox{or}}}(4,u_7)$,
$u_7\in\{4,5,6\}$.

\itemitem{(2)} $(a,b)=(3,2)$, $\{c,d\}=\{r,s\}$,
where $(r,s)=(3,u_8)$, $u_8\in\{3,4,5\}$.

\itemitem{(3)} $(a,b)=(2,r)$, $\{c,d\}=\{3,s\}$,
where $(r,s)=(4,u_9)$, $u_9\in\{4,5,\ldots,10\}$,
$\mathrel{\mathop{=\hskip-5pt=}\limits^{\hbox{or}}}(5,u_{10})$,
$u_{10}\in\{5,6,7,8\}$,
$\mathrel{\mathop{=\hskip-5pt=}\limits^{\hbox{or}}}(6,u_{11})$,
$u_{11}\in\{6,7,8\}$,
$\mathrel{\mathop{=\hskip-5pt=}\limits^{\hbox{or}}}(7,7)$.

\itemitem{(4)} $(a,b)=(2,s)$, $\{c,d\}=\{3,r\}$,
where $(s,r)=(s,3)$, $s\geq 4$,
$\mathrel{\mathop{=\hskip-5pt=}\limits^{\hbox{or}}}(s,4)$,
$s\geq 5$,
$\mathrel{\mathop{=\hskip-5pt=}\limits^{\hbox{or}}}(s,5)$,
$s\geq 6$,
$\mathrel{\mathop{=\hskip-5pt=}\limits^{\hbox{or}}}(s,6)$,
$s\geq 7$,
$\mathrel{\mathop{=\hskip-5pt=}\limits^{\hbox{or}}}(u_{12},7)$,
$u_{12}\in\{8,9,10\}$.

\itemitem{(5)} $(a,b)=(4,2)$, $\{c,d\}=\{3,3\}$.

\itemitem{(6)} $(a,b)=(3,r)$, $\{c,d\}=\{2,s\}$,
where $(r,s)=(3,u_{13})$, $u_{13}\in\{3,4,\ldots,13\}$,
$\mathrel{\mathop{=\hskip-5pt=}\limits^{\hbox{or}}}(4,u_{14})$,
$u_{14}\in\{4,5,\ldots,9\}$,
$\mathrel{\mathop{=\hskip-5pt=}\limits^{\hbox{or}}}(5,u_{15})$,
$u_{15}\in\{5,6,7,8\}$,
$\mathrel{\mathop{=\hskip-5pt=}\limits^{\hbox{or}}}(6,u_{16})$,
$u_{16}\in\{6,7\}$,
$\mathrel{\mathop{=\hskip-5pt=}\limits^{\hbox{or}}}(7,7)$.

\itemitem{(7)} $(a,b)=(r,3)$, $\{c,d\}=\{2,s\}$,
where $(r,s)=(4,u_{17})$, $u_{17}\in\{4,5\}$.

\itemitem{(8)} $(a,b)=(3,s)$, $\{c,d\}=\{2,r\}$,
where $(s,r)=(s,3)$, $s\geq 4$,
$\mathrel{\mathop{=\hskip-5pt=}\limits^{\hbox{or}}}(s,4)$,
$s\geq 5$,
$\mathrel{\mathop{=\hskip-5pt=}\limits^{\hbox{or}}}(s,5)$,
$s\geq 6$,
$\mathrel{\mathop{=\hskip-5pt=}\limits^{\hbox{or}}}(s,6)$,
$s\geq 7$,
$\mathrel{\mathop{=\hskip-5pt=}\limits^{\hbox{or}}}(u_{18},7)$,
$u_{18}\in\{8,9\}$.

\itemitem{(9)} $(a,b)=(s,3)$, $\{c,d\}=\{2,r\}$,
where $(s,r)=(u_{19},3)$, $u_{19}\in\{4,5,\ldots,8\}$,
$\mathrel{\mathop{=\hskip-5pt=}\limits^{\hbox{or}}}(5,4)$.

\itemitem{(10)} $(a,b)=(r,s)$, $\{c,d\}=\{2,3\}$,
where $(r,s)=(4,s)$, $s\geq 4$,
$\mathrel{\mathop{=\hskip-5pt=}\limits^{\hbox{or}}}(5,s)$,
$s\geq 5$,
$\mathrel{\mathop{=\hskip-5pt=}\limits^{\hbox{or}}}(6,s)$,
$s\geq 6$.

\itemitem{(11)} $(a,b)=(s,r)$, $\{c,d\}=\{2,3\}$,
where $(s,r)=(u_{20},4)$, $u_{20}\in\{5,6,7\}$,
$\mathrel{\mathop{=\hskip-5pt=}\limits^{\hbox{or}}}(u_{21},5)$,
$u_{21}\in\{6,7\}$.

\smallskip
\item{4.}
\Itemitem{(1)} $(a,b)=(2,4)$, $\{c,d\}=\{r,s\}$,
where $(r,s)=(4,u_{22})$, $u_{22}\in\{4,5\}$.

\itemitem{(2)} $(a,b)=(2,5)$, $\{c,d\}=\{4,5\}$.

\itemitem{(3)} $(a,b)=(2,s)$, $\{c,d\}=\{4,r\}$,
where $(s,r)=(s,4)$, $s\geq 5$.

\itemitem{(4)} $(a,b)=(4,r)$, $\{c,d\}=\{2,s\}$,
where $(r,s)=(4,u_{23})$, $u_{23}\in\{4,5\}$.

\itemitem{(5)} $(a,b)=(4,s)$, $\{c,d\}=\{2,r\}$,
where $(s,r)=(s,4)$, $s\geq 5$.

\smallskip
\item{5.}
\Itemitem{(1)} $(a,b)=(3,3)$, $\{c,d\}=\{r,s\}$,
where $(r,s)=(3,u_{24})$, $u_{24}\in\{3,4\}$.

\itemitem{(2)} $(a,b)=(3,s)$, $\{c,d\}=\{3,r\}$,
where $(s,r)=(s,3)$, $s\geq 4$.

\medskip

\noindent (XIV) $f(x,y,z,w)=x^a+xy^b+xz^c+xw^d+y^pz^q+z^rw^s$,
$\frac{p(a-1)}{ab}+\frac{q(a-1)}{ac}=1$,

\hskip10pt
$\frac{r(a-1)}{ac}+\frac{s(a-1)}{ad}=1$.

\smallskip
\item{1.}
\Itemitem{(1)}
$(a,b,c,d)=(r,1,s,t)$, $r\geq 2$, $s\geq 2$, $t\geq 2$.
\itemitem{(2)} $(a,b,c,d)=(r,s,1,t)$, $r\geq 2$, $s\geq 2$, $t\geq 2$.
\itemitem{(3)} $(a,b,c,d)=(r,s,t,1)$, $r\geq 2$, $s\geq 2$, $t\geq 2$.
\itemitem{(4)} $(a,b,c,d)=(r,1,1,s)$, $r\geq 2$, $s\geq 2$.
\itemitem{(5)} $(a,b,c,d)=(r,1,s,1)$, $r\geq 2$, $s\geq 2$.
\itemitem{(6)} $(a,b,c,d)=(r,s,1,1)$, $r\geq 2$, $s\geq 2$.
\itemitem{(7)} $(a,b,c,d)=(r,1,1,1)$, $r\geq 2$.

\smallskip
\item{2.}
\Itemitem{(1)} $a=2$, $\{b,c,d\}=\{2,r,s\}$,
where $(r,s)=(2,s)$, $s\geq 2$,
$\mathrel{\mathop{=\hskip-5pt=}\limits^{\hbox{or}}}(3,u_1)$,
$u_1\in\{3,4,5\}$.

\itemitem{(2)} $a=r$, $\{b,c,d\}=\{2,2,s\}$, $r\geq 3$, $s\geq r$.

\itemitem{(3)} $a=s$, $\{b,c,d\}=\{2,2,r\}$, $r\geq 2$, $s\geq r+1$.

\smallskip
\item{3.}
\Itemitem{(1)} $a=3$, $\{b,c,d\}=\{2,r,s\}$,
where $(r,s)=(3,u_2)$, $u_2\in\{3,4,5\}$.

\itemitem{(2)} $a=r$, $\{b,c,d\}=\{2,3,s\}$,
where $(r,s)=(4,u_3)$, $u_3\in\{4,5\}$,
$\mathrel{\mathop{=\hskip-5pt=}\limits^{\hbox{or}}}(5,5)$.

\itemitem{(3)} $a=s$, $\{b,c,d\}=\{2,3,r\}$,
where $(s,r)=(s,3)$, $s\geq 4$,
$\mathrel{\mathop{=\hskip-5pt=}\limits^{\hbox{or}}}(s,4)$,
$s\geq 5$,
$\mathrel{\mathop{=\hskip-5pt=}\limits^{\hbox{or}}}(s,5)$,
$y\geq 6$.

\medskip

\noindent (XV) $f(x,y,z,w)=x^ay+xy^b+xz^c+zw^d+y^pz^q$,
$\frac{p(a-1)}{ab-1}+\frac{qb(a-1)}{c(ab-1)}=1$.

\smallskip
\item{1.}
\Itemitem{(1)}
$(a,b,c,d)=(r,s,1,t)$, $r\geq 2$, $s\geq 2$, $t\geq 2$.
\itemitem{(2)} $(a,b,c,d)=(r,s,t,1)$, $r\geq 2$, $s\geq 2$, $t\geq 2$.
\itemitem{(3)} $(a,b,c,d)=(r,s,1,1)$, $r\geq 2$, $s\geq 2$.

\smallskip
\item{2.}
\Itemitem{(1)} $(a,b,c,d)=(2,2,r,s)$,
where $(r,s)=(2,s)$, $s\geq 2$,
$\mathrel{\mathop{=\hskip-5pt=}\limits^{\hbox{or}}}(3,u_1)$,
$u_1\in\{3,4,5,6\}$,
$\mathrel{\mathop{=\hskip-5pt=}\limits^{\hbox{or}}}(4,4)$.

\itemitem{(2)} $(a,b,c,d)=(2,2,s,r)$,
where $(s,r)=(s,2)$, $s\geq 3$,
$\mathrel{\mathop{=\hskip-5pt=}\limits^{\hbox{or}}}(s,3)$,
$s\geq 4$,
$\mathrel{\mathop{=\hskip-5pt=}\limits^{\hbox{or}}}(5,4)$.

\itemitem{(3)} $(a,b,c,d)=(2,r,2,s)$,
where $(r,s)=(3,u_2)$, $u_2\in\{3,4,5,6\}$,
$\mathrel{\mathop{=\hskip-5pt=}\limits^{\hbox{or}}}(4,4)$.

\itemitem{(4)} $(a,b,c,d)=(2,r,s,2)$, $r\geq 3$, $s\geq r$.

\itemitem{(5)} $(a,b)=(r,2)$, $\{c,d\}=\{2,s\}$, $r\geq 3$, $s\geq r$.

\itemitem{(6)} $(a,b,c,d)=(2,s,2,r)$,
where $(s,r)=(s,2)$, $s\geq 3$,
$\mathrel{\mathop{=\hskip-5pt=}\limits^{\hbox{or}}}(s,3)$,
$s\geq 4$,
$\mathrel{\mathop{=\hskip-5pt=}\limits^{\hbox{or}}}(5,4)$.

\itemitem{(7)} $(a,b,c,d)=(2,s,r,2)$, $r\geq 3$, $s\geq r+1$.

\itemitem{(8)} $(a,b)=(s,2)$, $\{c,d\}=\{2,r\}$, $r\geq 2$, $s\geq r+1$.

\itemitem{(9)} $(a,b,c,d)=(r,s,2,2)$,
where $(r,s)=(3,s)$, $s\geq 3$,
$\mathrel{\mathop{=\hskip-5pt=}\limits^{\hbox{or}}}(4,u_3)$,
$u_3\in\{4,5,\ldots,9\}$,
$\mathrel{\mathop{=\hskip-5pt=}\limits^{\hbox{or}}}(5,u_4)$,
$u_4\in\{5,6\}$.

\itemitem{(10)} $(a,b,c,d)=(s,r,2,2)$,
where $(s,r)=(s,3)$, $s\geq 4$,
$\mathrel{\mathop{=\hskip-5pt=}\limits^{\hbox{or}}}(s,4)$,
$s\geq 5$,
$\mathrel{\mathop{=\hskip-5pt=}\limits^{\hbox{or}}}(u_5,5)$,
$u_5\in\{6,7,8\}$.

\smallskip
\item{3.}
\Itemitem{(1)} $(a,b,c,d)=(2,3,3,3)$.

\itemitem{(2)} $(a,b,c,d)=(2,3,s,r)$,
where $(s,r)=(u_6,3)$, $u_6\in\{4,5\}$.

\itemitem{(3)} $(a,b,c,d)=(3,2,r,s)$,
where $(r,s)=(3,u_7)$, $u_7\in\{3,4,5\}$.

\itemitem{(4)} $(a,b,c,d)=(3,2,s,r)$,
where $(s,r)=(u_8,3)$, $u_8\in\{4,5,6,7\}$.

\itemitem{(5)} $(a,b,c,d)=(4,2,3,4)$.

\itemitem{(6)} $(a,b,c,d)=(r,2,s,3)$,
where $(r,s)=(4,u_9)$, $u_9\in\{4,5\}$,
$\mathrel{\mathop{=\hskip-5pt=}\limits^{\hbox{or}}}(5,5)$.

\itemitem{(7)} $(a,b,c,d)=(2,s,3,r)$,
where $(s,r)=(u_{10},3)$, $u_{10}\in\{4,5\}$.

\itemitem{(8)} $(a,b,c,d)=(s,2,3,r)$,
where $(s,r)=(s,3)$, $s\geq 4$,
$\mathrel{\mathop{=\hskip-5pt=}\limits^{\hbox{or}}}(s,4)$,
$s\geq 5$.

\itemitem{(9)} $(a,b,c,d)=(s,2,r,3)$,
where $(s,r)=(s,4)$, $s\geq 5$.

\itemitem{(10)} $(a,b,c,d)=(3,r,2,s)$,
where $(r,s)=(3,u_{11})$, $u_{11}\in\{3,4\}$.

\itemitem{(11)} $(a,b,c,d)=(3,r,s,2)$,
where $(r,s)=(3,s)$, $s\geq 3$,
$\mathrel{\mathop{=\hskip-5pt=}\limits^{\hbox{or}}}(4,u_{12})$,
$u_{12}\in\{4,5,6,7\}$.

\itemitem{(12)} $(a,b,c,d)=(4,3,2,4)$.

\itemitem{(13)} $(a,b,c,d)=(r,3,s,2)$,
where $(r,s)=(4,u_{13})$, $u_{13}\in\{4,5,\ldots,8\}$,
$\mathrel{\mathop{=\hskip-5pt=}\limits^{\hbox{or}}}(5,5)$.

\itemitem{(14)} $(a,b,c,d)=(3,4,2,3)$.

\itemitem{(15)} $(a,b,c,d)=(3,s,r,2)$,
where $(s,r)=(u_{14},3)$, $u_{14}\in\{4,5,\ldots,8\}$,
$\mathrel{\mathop{=\hskip-5pt=}\limits^{\hbox{or}}}(5,4)$.

\itemitem{(16)} $(a,b,c,d)=(s,3,2,r)$,
where $(s,r)=(s,3)$, $s\geq 4$.

\itemitem{(17)} $(a,b,c,d)=(s,3,r,2)$,
where $(s,r)=(s,3)$, $s\geq 4$,
$\mathrel{\mathop{=\hskip-5pt=}\limits^{\hbox{or}}}(u_{15},4)$,
$u_{15}\in\{5,6,7,8\}$.

\itemitem{(18)} $(a,b,c,d)=(4,4,3,2)$.

\itemitem{(19)} $(a,b,c,d)=(5,4,3,2)$.

\medskip

\noindent (XVI) $f(x,y,z,w)=x^ay+xy^b+xz^c+xw^d+y^pz^q+z^rw^s$,
$\frac{p(a-1)}{ab-1}+\frac{qb(a-1)}{c(ab-1)}=1$,

\hskip10pt $\frac{rb(a-1)}{c(ab-1)}+\frac{sb(a-1)}{d(ab-1)}=1$.

\smallskip

\item{1.}
\Itemitem{(1)}
$(a,b,c,d)=(r,s,1,t)$, $r\geq 2$, $s\geq 2$, $t\geq 2$.
\itemitem{(2)} $(a,b,c,d)=(r,s,t,1)$, $r\geq 2$, $s\geq 2$, $t\geq 2$.

\itemitem{(3)} $(a,b,c,d)=(r,s,1,1)$, $r\geq 2$, $s\geq 2$.

\smallskip
\item{2.}
\Itemitem{(1)} $\{a,b\}=\{2,2\}$, $\{c,d\}=\{r,s\}$,
where $(r,s)=(2,s)$, $s\geq 2$,
$\mathrel{\mathop{=\hskip-5pt=}\limits^{\hbox{or}}}(3,u_1)$,
$u_1\in\{3,4,5\}$.

\itemitem{(2)} $(a,b)=(2,r)$, $\{c,d\}=\{2,s\}$,
where $(r,s)=(3,u_2)$, $u_2\in\{3,4,5\}$.

\itemitem{(3)} $(a,b)=(r,2)$, $\{c,d\}=\{2,s\}$, $r\geq 3$, $s\geq r$.

\itemitem{(4)} $(a,b)=(2,s)$, $\{c,d\}=\{2,r\}$,
where $(s,r)=(s,2)$, $s\geq 3$,
$\mathrel{\mathop{=\hskip-5pt=}\limits^{\hbox{or}}}(u_3,3)$,
$u_3\in\{4,5\}$.

\itemitem{(5)} $(a,b)=(s,2)$, $\{c,d\}=\{2,r\}$, $r\geq 2$, $s\geq r+1$.

\itemitem{(6)} $\{a,b\}=\{r,s\}$, $\{c,d\}=\{2,2\}$, $r\geq 3$, $s\geq r$.

\smallskip
\item{3.}
\Itemitem{(1)} $(a,b)=(3,2)$, $\{c,d\}=\{r,s\}$,
where $(r,s)=(3,u_4)$, $u_4\in\{3,4,5\}$.

\itemitem{(2)} $(a,b)=(r,2)$, $\{c,d\}=\{3,s\}$,
where $(r,s)=(4,u_5)$, $u_5\in\{4,5\}$,
$\mathrel{\mathop{=\hskip-5pt=}\limits^{\hbox{or}}}(5,5)$.

\itemitem{(3)} $(a,b)=(s,2)$, $\{c,d\}=\{3,r\}$,
where $(s,r)=(s,3)$, $s\geq 4$,
$\mathrel{\mathop{=\hskip-5pt=}\limits^{\hbox{or}}}(s,4)$,
$s\geq 5$,
$\mathrel{\mathop{=\hskip-5pt=}\limits^{\hbox{or}}}(s,5)$,
$s\geq 6$.

\itemitem{(4)} $(a,b)=(3,r)$, $\{c,d\}=\{2,s\}$,
where $(r,s)=(3,u_6)$, $u_6\in\{3,4,5\}$.

\itemitem{(5)} $(a,b)=(r,3)$, $\{c,d\}=\{2,s\}$,
where $(r,s)=(4,u_7)$, $u_7\in\{4,5\}$,
$\mathrel{\mathop{=\hskip-5pt=}\limits^{\hbox{or}}}(5,5)$.

\itemitem{(6)} $(a,b)=(3,s)$, $\{c,d\}=\{2,r\}$,
where $(s,r)=(u_8,3)$, $u_8\in\{4,5\}$.

\itemitem{(7)} $(a,b)=(s,3)$, $\{c,d\}=\{2,r\}$,
where $(s,r)=(s,3)$, $s\geq 4$,
$\mathrel{\mathop{=\hskip-5pt=}\limits^{\hbox{or}}}(s,4)$,
$s\geq 5$,
$\mathrel{\mathop{=\hskip-5pt=}\limits^{\hbox{or}}}(s,5)$,
$s\geq 6$.

\itemitem{(8)} $(a,b)=(r,s)$, $\{c,d\}=\{2,3\}$,
where $(r,s)=(4,u_9)$, $u_9\in\{4,5\}$,
$\mathrel{\mathop{=\hskip-5pt=}\limits^{\hbox{or}}}(5,5)$.

\itemitem{(9)} $(a,b)=(s,r)$, $\{c,d\}=\{2,3\}$,
where $(s,r)=(s,4)$, $s\geq 5$,
$\mathrel{\mathop{=\hskip-5pt=}\limits^{\hbox{or}}}(s,5)$,
$y\geq 6$.

\medskip

\noindent (XVII) $f(x,y,z,w)=x^ay+xy^b+yz^c+xw^d+y^pz^q+z^rw^s$,
$\frac{p(a-1)}{ab-1}+\frac{qb(a-1)}{d(ab-1)}=1$,

\hskip10pt $\frac{r(b-1)}{ab-1}+\frac{sa(b-1)}{c(ab-1)}=1$.

\smallskip

\item{1.}
\Itemitem{(1)}
$(a,b,c,d)=(r,s,1,t)$, $r\geq 2$, $s\geq 2$, $t\geq 2$.
\itemitem{(2)} $(a,b,c,d)=(r,s,t,1)$, $r\geq 2$, $s\geq 2$, $t\geq 2$.
\itemitem{(3)} $(a,b,c,d)=(r,s,1,1)$, $r\geq 2$, $s\geq 2$.

\smallskip
\item{2.}
\Itemitem{(1)} $(a,b,c,d)=(2,2,r,s)$,
where $(r,s)=(2,s)$, $s\geq 2$,
$\mathrel{\mathop{=\hskip-5pt=}\limits^{\hbox{or}}}(3,u_1)$,
$u_1\in\{3,4,5\}$.

\itemitem{(2)} $(a,b,c,d)=(2,2,s,r)$,
where $(s,r)=(s,2)$, $s\geq 3$,
$\mathrel{\mathop{=\hskip-5pt=}\limits^{\hbox{or}}}(u_2,3)$,
$u_2\in\{4,5\}$.

\itemitem{(3)} $(a,b,c,d)=(2,r,2,s)$, $r\geq 3$, $s\geq r$.

\itemitem{(4)} $(a,b,c,d)=(2,r,s,2)$,
where $(r,s)=(3,u_3)$, $u_3\in\{3,4,\ldots,7\}$,
$\mathrel{\mathop{=\hskip-5pt=}\limits^{\hbox{or}}}(4,u_3)$,
$u_4\in\{4,5\}$,
$\mathrel{\mathop{=\hskip-5pt=}\limits^{\hbox{or}}}(5,5)$.

\itemitem{(5)} $(a,b,c,d)=(r,2,2,s)$,
where $(r,s)=(3,u_5)$, $u_5\in\{3,4,\ldots,7\}$,
$\mathrel{\mathop{=\hskip-5pt=}\limits^{\hbox{or}}}(4,u_6)$,
$u_6\in\{4,5\}$,
$\mathrel{\mathop{=\hskip-5pt=}\limits^{\hbox{or}}}(5,5)$.

\itemitem{(6)} $(a,b,c,d)=(r,2,s,2)$, $r\geq 3$, $s\geq r$.

\itemitem{(7)} $(a,b,c,d)=(2,s,2,r)$, $r\geq 2$, $s\geq r+1$.

\itemitem{(8)} $(a,b,c,d)=(2,s,r,2)$,
where $(s,r)=(s,3)$, $s\geq 4$,
$\mathrel{\mathop{=\hskip-5pt=}\limits^{\hbox{or}}}(s,4)$,
$s\geq 5$.

\itemitem{(9)} $(a,b,c,d)=(s,2,2,r)$,
where $(s,r)=(s,2)$, $s\geq 3$,
$\mathrel{\mathop{=\hskip-5pt=}\limits^{\hbox{or}}}(s,3)$,
$s\geq 4$,
$\mathrel{\mathop{=\hskip-5pt=}\limits^{\hbox{or}}}(s,4)$,
$s\geq 5$.

\itemitem{(10)} $(a,b,c,d)=(s,2,r,2)$, $r\geq 3$, $s\geq r+1$.

\itemitem{(11)} $\{a,b\}=\{r,s\}$, $\{c,d\}=\{2,2\}$, $r\geq 3$, $s\geq r$.

\smallskip
\item{3.}
\Itemitem{(1)} $(a,b,c,d)=(2,3,r,s)$,
where $(r,s)=(3,u_7)$, $u_7\in\{3,4\}$.

\itemitem{(2)} $(a,b,c,d)=(3,2,3,3)$.

\itemitem{(3)} $(a,b,c,d)=(3,2,4,3)$.

\itemitem{(4)} $(a,b,c,d)=(2,s,3,r)$,
where $(s,r)=(s,3)$, $s\geq 4$.

\itemitem{(5)} $(a,b,c,d)=(s,2,3,r)$,
where $(s,r)=(s,3)$, $s\geq 4$.

\itemitem{(6)} $(a,b,c,d)=(3,r,2,s)$,
where $(r,s)=(3,u_8)$, $u_8\in\{3,4,5\}$,
$\mathrel{\mathop{=\hskip-5pt=}\limits^{\hbox{or}}}(4,u_{9})$,
$u_9\in\{4,5\}$.

\itemitem{(7)} $(a,b,c,d)=(3,r,s,2)$,
where $(r,s)=(3,u_{10})$, $u_{10}\in\{3,4,5\}$,
$\mathrel{\mathop{=\hskip-5pt=}\limits^{\hbox{or}}}(4,4)$.

\itemitem{(8)} $(a,b,c,d)=(4,3,2,4)$.

\itemitem{(9)} $(a,b,c,d)=(r,3,s,2)$,
where $(r,s)=(4,u_{11})$, $u_{11}\in\{4,5\}$.

\itemitem{(10)} $(a,b,c,d)=(3,s,2,r)$,
where $(s,r)=(s,3)$, $s\geq 4$,
$\mathrel{\mathop{=\hskip-5pt=}\limits^{\hbox{or}}}(s,4)$,
$s\geq 5$.

\itemitem{(11)} $(a,b,c,d)=(3,s,r,2)$,
where $(s,r)=(s,3)$, $s\geq 4$.

\itemitem{(12)} $(a,b,c,d)=(s,3,2,r)$,
where $(s,r)=(s,3)$, $s\geq 4$.

\itemitem{(13)} $(a,b,c,d)=(s,3,r,2)$,
where $(s,r)=(s,3)$, $s\geq 4$,
$\mathrel{\mathop{=\hskip-5pt=}\limits^{\hbox{or}}}(s,4)$,
$s\geq 5$.

\itemitem{(14)} $(a,b,c,d)=(r,s,2,3)$,
where $(r,s)=(4,s)$, $s\geq 4$,
$\mathrel{\mathop{=\hskip-5pt=}\limits^{\hbox{or}}}(5,u_{12})$,
$u_{12}\in\{5,6,7,8\}$.

\itemitem{(15)} $(a,b,c,d)=(r,s,3,2)$,
where $(r,s)=(4,u_{13})$, $u_{13}\in\{4,5,\ldots,9\}$,
$\mathrel{\mathop{=\hskip-5pt=}\limits^{\hbox{or}}}(5,u_{14})$,
$u_{14}\in\{5,6\}$.

\itemitem{(16)} $(a,b,c,d)=(s,r,2,3)$,
where $(s,r)=(u_{15},4)$, $u_{15}\in\{5,6,\ldots,9\}$,
$\mathrel{\mathop{=\hskip-5pt=}\limits^{\hbox{or}}}(6,5)$.

\itemitem{(17)} $(a,b,c,d)=(s,r,3,2)$,
where $(s,r)=(s,4)$, $s\geq 5$,
$\mathrel{\mathop{=\hskip-5pt=}\limits^{\hbox{or}}}(u_{16},5)$,
$u_{16}\in\{6,7,8\}$.

\medskip

\noindent (XVIII) $f(x,y,z,w)=x^az+xy^b+yz^c+yw^d+z^pw^q$,
$\frac{p(a(b-1)+1)}{abc+1}+\frac{qc(a(b-1)+1)}{d(abc+1)}=1$.

\smallskip
\item{1.}
\Itemitem{(1)}
$(a,b,c,d)=(1,r,s,t)$, $r\geq 2$, $s\geq 2$, $t\geq 2$.
\itemitem{(2)} $(a,b,c,d)=(r,1,s,t)$, $r\geq 2$, $s\geq 2$, $t\geq 2$.

\itemitem{(3)} $(a,b,c,d)=(r,s,1,t)$, $r\geq 2$, $s\geq 2$, $t\geq 2$.
\itemitem{(4)} $(a,b,c,d)=(r,s,t,1)$, $r\geq 2$, $s\geq 2$, $t\geq 2$.
\itemitem{(5)} $(a,b,c,d)=(1,1,r,s)$, $r\geq 2$, $s\geq 2$.
\itemitem{(6)} $(a,b,c,d)=(1,r,1,s)$, $r\geq 2$, $s\geq 2$.
\itemitem{(7)} $(a,b,c,d)=(1,r,s,1)$, $r\geq 2$, $s\geq 2$.
\itemitem{(8)} $(a,b,c,d)=(r,1,1,s)$, $r\geq 2$, $s\geq 2$.
\itemitem{(9)} $(a,b,c,d)=(r,1,s,1)$, $r\geq 2$, $s\geq 2$.
\itemitem{(10)} $(a,b,c,d)=(r,s,1,1)$, $r\geq 2$, $s\geq 2$.
\itemitem{(11)} $(a,b,c,d)=(1,1,1,r)$, $r\geq 2$.
\itemitem{(12)} $(a,b,c,d)=(1,1,r,1)$, $r\geq 2$.
\itemitem{(13)} $(a,b,c,d)=(1,r,1,1)$, $r\geq 2$.
\itemitem{(14)} $(a,b,c,d)=(r,1,1,1)$, $r\geq 2$.
\itemitem{(15)} $(a,b,c,d)=(1,1,1,1)$.

\smallskip
\item{2.}
\Itemitem{(1)} $(a,b,c,d)=(2,2,r,s)$,
where $(r,s)=(2,s)$, $s\geq 2$,
$\mathrel{\mathop{=\hskip-5pt=}\limits^{\hbox{or}}}(3,u_1)$,
$u_1\in\{3,4,\ldots,8\}$,
$\mathrel{\mathop{=\hskip-5pt=}\limits^{\hbox{or}}}(4,u_2)$,
$u_2\in\{4,5\}$.

\itemitem{(2)} $(a,b,c,d)=(2,2,s,r)$,
where $(s,r)=(s,2)$, $s\geq 3$,
$\mathrel{\mathop{=\hskip-5pt=}\limits^{\hbox{or}}}(s,3)$,
$s\geq 4$,
$\mathrel{\mathop{=\hskip-5pt=}\limits^{\hbox{or}}}(u_3,4)$,
$u_4\in\{5,6,7\}$.

\itemitem{(3)} $(a,b,c,d)=(2,r,2,s)$,
where $(r,s)=(3,u_4)$, $u_4\in\{3,4,\ldots,9\}$,
$\mathrel{\mathop{=\hskip-5pt=}\limits^{\hbox{or}}}(4,u_5)$,
$u_5\in\{4,5,6\}$,
$\mathrel{\mathop{=\hskip-5pt=}\limits^{\hbox{or}}}(5,5)$.

\itemitem{(4)} $(a,b,c,d)=(2,r,s,2)$, $r\geq 3$, $s\geq r$.

\itemitem{(5)} $(a,b,c,d)=(r,2,2,s)$,
where $(r,s)=(3,u_6)$, $u_6\in\{3,4,5,6,7\}$,
$\mathrel{\mathop{=\hskip-5pt=}\limits^{\hbox{or}}}(4,4)$.

\itemitem{(6)} $(a,b,c,d)=(r,2,s,2)$,
where $(r,s)=(3,s)$, $s\geq 3$,
$\mathrel{\mathop{=\hskip-5pt=}\limits^{\hbox{or}}}(4,u_7)$,
$u_7\in\{4,5,6,7\}$.
\itemitem{(7)} $(a,b,c,d)=(2,s,2,r)$,
where $(s,r)=(s,2)$, $s\geq 3$,
$\mathrel{\mathop{=\hskip-5pt=}\limits^{\hbox{or}}}(s,3)$,
$s\geq 4$,
$\mathrel{\mathop{=\hskip-5pt=}\limits^{\hbox{or}}}(s,4)$,
$s\geq 5$,
$\mathrel{\mathop{=\hskip-5pt=}\limits^{\hbox{or}}}(u_8,5)$,
$u_8\in\{6,7\}$.

\itemitem{(8)} $(a,b,c,d)=(2,s,r,2)$, $r\geq 3$, $s\geq r+1$.

\itemitem{(9)} $(a,b,c,d)=(s,2,2,r)$,
where $(s,r)=(s,2)$, $s\geq 3$,
$\mathrel{\mathop{=\hskip-5pt=}\limits^{\hbox{or}}}(u_9,3)$,
$u_9\in\{4,5,6,7\}$.

\itemitem{(10)} $(a,b,c,d)=(s,2,r,2)$,
where $(s,r)=(u_{10},3)$, $u_{10}\in\{4,5,\ldots,8\}$,
$\mathrel{\mathop{=\hskip-5pt=}\limits^{\hbox{or}}}(5,4)$.

\itemitem{(11)} $\{a,b\}=\{r,s\}$, $\{c,d\}=\{2,2\}$, $r\geq 3$, $s\geq r$.

\smallskip
\item{3.}
\Itemitem{(1)} $(a,b,c,d)=(2,3,r,s)$,
where $(r,s)=(3,u_{11})$, $u_{11}\in\{3,4\}$.

\itemitem{(2)} $(a,b,c,d)=(2,3,s,r)$,
where $(s,r)=(u_{12},3)$, $u_{12}\in\{4,5,\ldots,8\}$.

\itemitem{(3)} $(a,b,c,d)=(3,2,3,3)$.

\itemitem{(4)} $(a,b,c,d)=(3,2,4,3)$.

\itemitem{(5)} $(a,b,c,d)=(2,4,3,4)$.

\itemitem{(6)} $(a,b,c,d)=(2,r,s,3)$,
where $(r,s)=(4,u_{13})$, $u_{13}\in\{4,5\}$.

\itemitem{(7)} $(a,b,c,d)=(2,s,3,r)$,
where $(s,r)=(s,3)$, $s\geq 4$.

\itemitem{(8)} $(a,b,c,d)=(2,s,r,3)$,
where $(s,r)=(u_{14},4)$, $u_{14}\in\{5,6,7\}$.

\itemitem{(9)} $(a,b,c,d)=(3,r,2,s)$,
where $(r,s)=(3,u_{15})$, $u_{15}\in\{3,4\}$.

\itemitem{(10)} $(a,b,c,d)=(3,r,s,2)$,
where $(r,s)=(3,u_{16})$, $u_{16}\in\{3,4,\ldots,9\}$,
$\mathrel{\mathop{=\hskip-5pt=}\limits^{\hbox{or}}}(4,u_{17})$,
$u_{17}\in\{4,5,6\}$,
$\mathrel{\mathop{=\hskip-5pt=}\limits^{\hbox{or}}}(5,5)$.

\itemitem{(11)} $(a,b,c,d)=(4,3,4,2)$.

\itemitem{(12)} $(a,b,c,d)=(3,s,2,r)$,
where $(s,r)=(s,3)$, $s\geq 4$.

\itemitem{(13)} $(a,b,c,d)=(3,s,r,2)$,
where $(s,r)=(s,3)$, $s\geq 4$,
$\mathrel{\mathop{=\hskip-5pt=}\limits^{\hbox{or}}}(s,4)$,
$s\geq 5$,
$\mathrel{\mathop{=\hskip-5pt=}\limits^{\hbox{or}}}(u_{18},5)$,
$u_{18}\in\{6,7\}$.

\itemitem{(14)} $(a,b,c,d)=(s,3,2,r)$,
where $(s,r)=(u_{19},3)$, $u_{19}\in\{4,5\}$.

\itemitem{(15)} $(a,b,c,d)=(s,3,r,2)$,
where $(s,r)=(u_{20},3)$, $u_{20}\in\{4,5,6\}$.

\itemitem{(16)} $(a,b,c,d)=(r,s,2,3)$,
where $(r,s)=(4,u_{21})$, $u_{21}\in\{4,5\}$.

\itemitem{(17)} $(a,b,c,d)=(r,s,3,2)$,
where $(r,s)=(4,s)$, $s\geq 4$,
$\mathrel{\mathop{=\hskip-5pt=}\limits^{\hbox{or}}}(5,5)$.

\itemitem{(18)} $(a,b,c,d)=(5,4,3,2)$.

\medskip

\noindent (XIX) $f(x,y,z,w)=x^az+xy^b+z^cw+yw^d$.

\smallskip

\item{1.}
\Itemitem{(1)}
$(a,b,c,d)=(1,r,s,t)$, $r\geq 2$, $s\geq 2$, $t\geq 2$.
\itemitem{(2)} $(a,b,c,d)=(1,1,r,s)$, $r\geq 2$, $s\geq 2$.

\item{2.}
\Itemitem{(1)} $(a,d)=(2,2)$, $(b,c)=(r,s)$, $r\geq 2$, $s\geq 2$.

\itemitem{(2)} $(a,d)=(2,r)$, $\{b,c\}=\{2,s\}$,
where $(r,s)=(3,s)$, $s\geq 3$,
$\mathrel{\mathop{=\hskip-5pt=}\limits^{\hbox{or}}}(4,u_1)$,
$u_1\in\{4,5,\ldots,8\}$.

\itemitem{(3)} $(a,d)=(2,5)$, $(b,c)=(2,5)$.

\smallskip
\item{3.}
\Itemitem{(1)} $(a,d)=(2,3)$, $\{b,c\}=\{r,s\}$,
where $(r,s)=(3,u_2)$, $u_2\in\{3,4,5,6,7\}$.

\itemitem{(2)} $(a,d)=(2,3)$, $(b,c)=(4,4)$.

\itemitem{(3)} $(a,d)=(2,s)$, $(b,c)=(3,r)$,
where $(s,r)=(u_3,3)$, $u_3\in\{4,5\}$.

\medskip

\noindent {\bf Proof of Theorem 3.3.}\quad
In view of Corollary 2.5 and Theorem 3.2, it is clear that an isolated
rational hypersurface singularity with $\c^*$-action is defined by one
of the 19 types in section 2 with $p_g=0$.
The equations of the $\Gamma_-$ hyperplanes of these 19 types are
respectively given by $\alpha(x,y,z,w)=1$.

\smallskip

In order to find all hypersurfaces among these 19 types with $p_g=0$,
we only need to find all solutions of $\alpha(1,1,1,1)>1$ among these
19 types.
We have used the MAPLE program [Ch] to perform the computations.
The solutions are listed in the statement of this Theorem.

\medskip

\noindent {\bf Remark.}\quad
The lists in (VIII), (XIII), (XIV), (XVI), (XVII) in Theorem 3.3 may be
reduced slightly by change of coordinates.

\bigskip
\bigskip

\centerline{\bf References}
\medskip
\noindent [Ar]\quad Arnold, V.I., Normal forms of functions in neighborhoods
of degenerate critical points, Russian

\hskip7pt Math. Surveys 29(1975), 10--50.

\smallskip
\noindent [Art]\quad Artin, M., On isolated rational singularities
of surfaces, Amer. J. Math. {\bf 88}(1966), 129--136.

\smallskip
\noindent [Bu]\quad Burns, D., On rational singularities in dimension $>2$,
Math. Ann., 211(1974), 237--244.

\smallskip
\noindent [Ch]\quad Char, B.W. [et al.] Maple V Language reference manual,
Springer Verlag, 1991.
\smallskip

\noindent [Ka]\quad Kannowski, M.A., Simply connected four manifolds
obtained from weighted homogeneous

\hskip8pt polynomials, Dissertation, The University
of Iowa, 1986.

\smallskip
\noindent [Ko]\quad Kouchnirenko, A.G., Poly\'{e}dres de Newton et Nombres
de Milnor, Invent. Math., 32(1976),

\hskip8pt 1--31.

\smallskip
\noindent [L-Y-Y]\quad Luk, H.S., Yau, S. S.-T., and Yu, Y.,
Algebraic classification and obstructions to

\hskip24pt embedding of strongly
pseudoconvex compact 3-dimensional CR manifolds in $\c^3$, Math.

\hskip24pt Nachr. 170(1994), 183--200.

\smallskip
\noindent [Me-Te]\quad Merle, M., and Teissier, B., Conditions d'Adjonction
d'apr\'{e}s Du Val, S\'{e}minaire sur les

\hskip24pt Singularit\'{e}s des Surfaces
(Centre de Math. de l'Ecole Polytechnique, 1976-1977),
Lecture

\hskip24pt Notes in Math., Vol. 777, Springer, Berlin, 1980, 229--245.

\smallskip
\noindent [Or-Ra]\quad Orlik, P. and Randell, R., The monodromy of weighted
homogeneous singularities, Invent.

\hskip24pt Math. 39, (1977), 199--211.

\smallskip
\noindent [Or-Wa]\quad Orlik, P., and Wagreich, P., Isolated singularities
of algebraic surfaces with $\c^*$-action,

\hskip24pt Ann. of Math. 93, (1971), 205--228.

\smallskip
\noindent [Xu-Ya1]\quad Xu, Y.-J., and Yau S. S.-T., Classification of
topological types of isolated quasi-homogeneous

\hskip30pt two dimensional hypersurface
singularities, Manuscripta Math. 64, (1989), 445--469.

\smallskip
\noindent [Xu-Ya2]\quad Xu, Y.-J., and Yau, S. S.-T., Topological types
of seven classes of isolated singularities

\hskip30pt with $\c^*$-action, Rocky
Mountain Journal of Mathematics, Vol. 22, (1992), 1147--1215.

\smallskip
\noindent [Ya1]\quad Yau, S. S.-T., Topological types of isolated hypersurface
singularities, Contemp. Math.

\hskip14pt 101(1989), 303--321.

\smallskip
\noindent [Ya2]\quad Yau, S. S.-T., Two theorems on higher dimensional
singularities, Math. Ann. 231(1977),

\hskip14pt 55--59.

\smallskip
\noindent [Ya3]\quad Yau, S. S.-T., Sheaf cohomology on 1-convex manifolds,
Recent Developments in Several

\hskip14pt Complex Variables, Ann. of Math. Study,
100(1981), 429--452.

\smallskip
\noindent [Ya-Yu]\quad Yau, S. S.-T. and Yu, Y., Algebraic classification
of rational CR structures on topological

\hskip24pt 5-sphere with transversal
holomorphic $S^1$-action in $\c^4$, Math. Nachrichten
246-247(2002),

\hskip24pt 207-233.

\bye